\documentclass{amsart}

\usepackage{amsmath,amssymb,amsthm}
\usepackage{geometry}
\usepackage[utf8]{inputenc}
\usepackage{color,graphicx}
\usepackage{enumitem}
\usepackage{mathrsfs}
\usepackage{bbm}
\usepackage[capitalize]{cleveref}
\usepackage{stmaryrd}
\usepackage{appendix}
\usepackage[numbers]{natbib}

\newcommand{\scp}[2]{{\big\langle {#1}\, , \, {#2}\big\rangle}}

\newcommand{\lform}[2]{{\big( {#1}\, \big| \, {#2}\big)}}

\newcommand{\Diff}{{\mathit{Diff}\!}}
\newcommand{\Sym}{{\mathit{Sym}\!}}

\newcommand{\CS}{{\mathscr{S}}}
\newcommand{\mR}{{\mathbb R}}

\newcommand{\id}[1]{\mathrm{id}_{\mathbb R^{#1}}}
\newcommand{\Id}[1]{\mathrm{Id}_{\mathbb R^{#1}}}
\newcommand{\dd}{{\mathbbmss d}}

\newcommand{\prt}{\partial}

\renewcommand{\phi}{\varphi}

\renewcommand{\epsilon}{\varepsilon}

\let\on=\operatorname

\newcommand{\wto}{\rightharpoonup}

\def\Emb{\mathit{Emb}}
\def\Diff{\mathit{Diff}}

\def\trace{\on{trace}}

\crefname{condition}{condition}{conditions}

\graphicspath{{./}{./normalCortex/}}

\title{Shape spaces: From geometry to biological plausibility}
\author[N. Charon]{Nicolas Charon}
\address{N. Charon. Center for Imaging Science, Johns Hopkins University}
\email{ncharon1@jhu.edu}
\author[L. Younes]{Laurent Younes}
\address{L. Younes. Center for Imaging Science, Johns Hopkins University}
\email{laurent.younes@jhu.edu}
\thanks{Nicolas Charon is partially supported by NSF 1945224 and NSF 1953267.}
\thanks{Laurent Younes is partially supported by NIH U19AG033655, R01NS102670 and R01AG055121.}

\begin{document}
\bibliographystyle{plainnat}

\begin{abstract}
    This chapter reviews several Riemannian metrics and evolution equations in the context of diffeomorphic shape analysis. After a short review of of various approaches at building Riemannian spaces of shapes, with a special focus on the foundations of the large deformation diffeomorphic metric mapping algorithm, the attention is turned to elastic metrics, and to growth models that can be derived from it. In the latter context, a new class of metrics, involving the optimization of a growth tensor, is introduced and some of its properties are studied.  
\end{abstract}
\maketitle

\section{Introduction: Shape spaces}

Shape has long been an object of scientific study, especially in life sciences where it provided a primary element in the differentiation between species. It was---in complement to behavioral patterns---a central factor of the early justification of evolutionary theory, and was of course the main subject of D'Arcy Thompson seminal work ``On Growth and Forms'' \citep{thompson1917growth}.

The construction of mathematical models of shape spaces, however, was more recent, and started with David Kendall's landmark paper introducing a shape space as a particular Riemannian manifold \citep{kendall1984shape}, a construction motivated by the need to provide a formal mathematical framework for statistical analyses of shape datasets. In Kendall's model, shapes are represented as  ordered collections of distinct points with fixed cardinality. The manifold structure is obtained as a quotient space through the action of rotations, translations and scaling and the metric as the projection of the Euclidean metric  to this quotient space.  Kendall's shape space has since been used in a large variety of applications, with increasing numbers of available shape datasets and relevant associated statistical questions (see the recent edition of \citet{dryden2016statistical} for additional details and references). 

Kendall's shape space is however limited by the need to provide a consistent ordering (or labeling) of the points constituting the shape, and by the requirement that they form a finite set. Shape datasets are typically formed by unlabelled geometric objects, and using Kendall's shape space requires defining and indexing (often manually) collections of landmarks for each shape, resulting in an intensive and sometimes imperfectly specified  problem. Defining shape spaces whose elements are curves or surfaces requires however more advanced mathematical tools, notably from global analysis \citep{palais1968foundations}, and a recent description of various formulations of shape spaces in this general context can be found in \citet{bauer2014overview}. In spite of the additional mathematical technicality, the construction of these shape spaces follows the general principles leading to Kendall's space: first define a simple space of geometric objects as an open subset of a normed (or Fr\'echet) space, where the finite-dimensional space of ordered distinct points is replaced, e.g., by a space of immersions (or embeddings) from a fixed manifold $M$ (the parameter space) to $\mR^d$, the ambient space.  This space (and its norm) is then quotiented by group actions to which shapes must be invariant, bringing in, in addition to previous actions of translations, rotations and scaling, the infinite-dimensional group of reparametrizations, provided by diffeomorphisms of the parameter space. Another modification to the finite-dimensional framework is that the Euclidean metric, as the base norm, which was a  natural choice when working with finite sets of points, now needs to be replaced with some invariant Hilbert metric (if one wants a Riemannian structure at the end) on the space of immersions, for which there are many choices, including the  whole family of invariant Sobolev norms. The well-posedness of various concepts in the resulting shape space, such as the non-degeneracy of the metric or the existence of geodesics, indeed depends on this choice. A striking example is the fact that the Riemannian distance between any pair of shapes may trivialize to zero for certain metrics, as initially discovered in \citet{michor2005vanishing} in the case of curves and then extended to other shapes spaces \citep{bauer2020vanishing}.
   
From the whole variety of shape spaces that can be built following this construction, a small number actually leads to practical algorithms and numerical implementations, which is an essential requirement when the goal is to analyze shape datasets. For curves, an important example is associated with a class of first-order Sobolev metrics on the space of immersions. One can indeed show that, after quotienting out rotations, translations and/or scalings, the resulting Riemannian manifold is isomorphic to standard manifolds (such as the infinite dimensional sphere, and Stiefel or Grassman manifolds) on which geodesic and geodesic distances can be explicitly computed. For curves, the additional cost of adding reparametrization invariance remains manageable, using, e.g., dynamic programming methods. A first example of such metrics was provided in \citet{younes1996distance,younes1998computable} with further developments in \citet{younes2008metric}. A second example was then provided in \citet{klassen2004analysis} (see \citet{srivastava2016functional}), and the approach was later extended to a one-parameter family including these two examples in \citet{needham2020simplifying} and \citet{younes2019shapes}, chapter 12 (see also \citet{bauer2014constructing}).  

Shape spaces have also been built using a different angle, leveraging the action of the diffeomorphism group of $\mR^d$ on a shape space. Diffeomorphisms of their ambient space indeed act transitively on most shapes of interest assuming that one fixes their topology (taking an example, diffeomorphisms of $\mR^2$ can be used to transform any $C^1$ Jordan plane curve to any other). Using a metric on the diffeomorphism group with suitable properties, one can, given two shapes, compute the diffeomorphism closest to the identity that transforms the first shape into the other, and the distance between the identity and this optimal diffeomorphism also provides a distance between the considered shapes. (This construction will be described in detail in \cref{sec:diffeo}.) Formally, the considered shape space is provided by all diffeomorphic transformation of  a given template. This approach can be seen as an application of Grenander's metric pattern theory \citep{grenander1991shape,grenander1993general} and as a mathematical formulation of D'Arcy Thompson models \citep{thompson1917growth}. It was introduced for shape spaces of images synchronously in \citet{dupuis1998variation} (with a precursor in \citet{christensen1996deformable}) and \citet{trouve1995action,trouve1998diffeomorphism}, and for collections of labelled points in \citet{miller1999large}. This formulation, very flexible, has later been applied to various shape spaces, such as unlabelled point sets \citep{glaunes2004diffeomorphic},  curves and surfaces \citep{vaillant2005surface,glaunes2008large}, vector or tensor fields \citep{cao2005large,cao2006diffeomorphic}. The reader may also refer to the recent survey in \citet{bauer2019metric} that describes in details the two previous approaches in the case of curves and surfaces.

Note that the previous discussion does not include the many methods that provide shape features, i.e., finite or infinite-dimensional descriptors that can be attached to a given shape,  without necessarily providing them with a clear mathematical structure (such as that of a Riemannian manifold) which is one of the main concerns of the construction of shape spaces. Such methods were introduced in computer vision, medical imaging, biology and are too numerous to cite exhaustively in this chapter. Among the most important ones (a subjective statement), one can cite approaches using complex analysis and the (quasi-) conformal maps to represent surfaces  \citep{gu2004genus,gu2008computational,zeng2011registration,zeng2012computing,lui2014shape}, isometry-invariant descriptors based on distance maps or Laplace-Beltrami eigenvectors \citep{bronstein2008analysis,bronstein2008numerical,memoli2008gromov-hausdorff,ovsjanikov2010one,memoli2011gromov}, or the shape context \citep{belongie2002shape}. The rest of this chapter will however remain focused on shape space approaches.

The construction of shape spaces as described above is based on purely geometric aspects. No physical law, or biological mechanism is used to define the various components that constitute the shape space. This non-committal approach is indeed justified, as shapes spaces are designed as containers for families of shapes that are not related to each other by a natural process (e.g., there is no physical process by which a finch's beak can transform into the shape of another one). This fact provides the technical advantage that the construction of shape metrics is not constrained by the laws of nature and can therefore be selected so that they guarantee the existence, say, of geodesics, provide nicely behaved gradient flows, etc... This will be illustrated in  \cref{sec:diffeo}.

On the other hand, biological processes provide many examples in which shapes change with time, in a process that is constrained by well specified laws. The goal of this chapter is to describe a few among recent attempts at representing such processes as trajectories in the shape spaces above, which, after small modifications or regularization, will be associated with evolution dynamics that behave well enough to allow for long time analysis and optimal control formulations.

This chapter is organized as follows. \Cref{sec:diffeo} provides a summary of the construction of shape spaces through diffeomorphic action. \Cref{sec:hybrid} focuses on variations of this construction with metrics that are inspired from elastic materials. \Cref{sec:growth} introduces a few examples of growth models in the context of shape spaces.
For an extensive introduction to mathematical models of growth, the reader should refer to \citet{goriely2017mathematics}, which provides a splendid reference on the topic, and in particular on ``morphoelasticity.'' The  representation of shape growth described in \cref{sec:growth} will, however, deviate to some extent from that described in this reference, and more generally from the large literature exploring morphoelasticity, as models will be designed in the form of control systems, with a control equation interpreted as a differential equation in shape space and growth or atrophy directly associated with the control. 
%We will also describe some preliminary attempts at inverting this process, that is, solving an optimal control problem optimizing over a small number of free parameters. This issue which, as we will see, is quite challenging, is nonetheless important since understanding plausible sources of atrophy may, for example, provide valuable hypotheses on the pathogenesis of neurodegenerative diseases. 

\section{Shape spaces under diffeomorphic action}
\label{sec:diffeo}

This section provides a summary of the construction of shape spaces based on the principles of D'Arcy-Thompson's theory of transformations \citep{thompson1917growth} and Grenander's metric pattern theory \citep{grenander1993general}. The fundamental principles of the construction were laid in \citet{dupuis1998variation,trouve1995action,grenander1998computational} and the reader may refer to \citet{younes2019shapes,miller_hamiltonian_2015,bauer2019metric} for more recent accounts of the theory. 

Shapes are modeled as embeddings from a fixed Riemannian manifold $M$ into $\mathbb R^d$, and therefore have a with fixed topology (in practice, $d=2$ or 3). Typically, $M$ is a unit circle or sphere, or a template shape of which one is computing deformations. Denote by $\Emb^p(M)$ the set of such $C^p$ embeddings, or simply $\Emb$ when $p$ and $M$ are fixed. Each element $m \in \Emb$ provides a shape equipped with a parametrization. Objects of interest are shapes modulo parametrization (also called ``unparametrized shapes'') in which one identifies embeddings $m$ and $\tilde m$ when they are related with each other through a change of parametrization, i.e., $\tilde m = m \circ \rho$ where $\rho$ is a diffeomorphism of $M$. In other terms, the shape space is defined as the quotient space of $\Emb$ through the right action of the diffeomorphism group of $M$, and will be denoted as $\CS$. Elements of $\CS$ will be denoted as $[m]$, for the equivalence class of $m\in \Emb$.

Comparisons between shapes rely on the group of transformations acting on $\Emb$ or $\CS$, which are modeled as diffeomorphisms of $\mR^d$. Denote by $\Diff^p(\mR^d)$, or simply $\Diff^p$, the group of $C^p$ diffeomorphisms of $\mR^d$ and by $\Diff^p_0(\mR^d)$, or simply $\Diff^p_0$, the subgroup of diffeomorphisms that converge to the identity map, denoted $\id{d}$, at infinity (convergence being understood in the $C^p$ sense). If $\phi\in \Diff^p$ and $m\in \Emb$, $\phi\cdot m$ is simply $\phi\circ m$, and this action commutes with reparametrization, so that one can define $\phi\cdot [m] = [\phi\cdot m]$ without ambiguity.

To compare two embeddings (or their associated shapes) $m$ and $m'$, one considers the transformations $\phi\in \Diff^p_0$ that relate them, i.e., such that $m' = \phi\cdot m$. One considers that $m$ and $m'$ are similar if one can find some $\phi$ relating them that is close to $\id{d}$. This closeness is itself evaluated using a metric on $\Diff^p_0$, with a  construction described below.

To provide a Riemannian metric, one needs an inner-product norm that evaluates the velocity of time-dependent diffeomorphisms, or ``diffeomorphic motions,'' taking the form $(x \mapsto \prt_t \phi(t, x))$ where $\phi$ is a function of time and space such that $(x\mapsto \phi(t,x))$ is at all times an element of $\Diff_0^p$. For a given time t, $(x\mapsto \prt_t\phi(t,x))$ is a $C^p$ vector field on $\mR^d$ and one therefore needs to provide a norm over such vector fields. The norm of the velocity at time $t$ should in principle depend on the diffeomorphism at the same time, $(x\mapsto \phi(t,x))$, but, for reasons seen below, it will be desirable for this norm to satisfy the invariance property that, when writing 
\[
\phi(t+\delta t, x) = \phi(t,x) + \delta\phi(t,x) = (\id{d} + \delta\phi(t, \cdot) \circ \phi^{-1}(t, \cdot))\circ \phi(t.x),
\]
the cost associated with $\delta \phi$ is a fixed function of the deformation increment $\delta\phi \circ \phi^{-1}$. Passing to the limit, this means that the Riemannian norm of $(x\mapsto \prt_t\phi(t,x))$ at $(x\mapsto \phi(t,x))$ is equal to the norm of $(x\mapsto \prt_t\phi(t,x) \circ \phi^{-1}(t,x))$ at $\id{d}$. The vector field $v(t,x) = \prt_t\phi(t,x) \circ \phi^{-1}(t,x)$ is called the Eulerian velocity of the diffeomorphic motion $x\mapsto \phi(t,x)$, and the diffeomorphic motion is recovered from the Eulerian velocity by solving the ordinary differential equation 
\begin{equation}
\label{eq:ode}
    \prt_t \phi(t,x) = v(t, \phi(t,x)).
\end{equation}

To define our Riemannian metric on $\Diff^p_0$, it therefore suffices to specify a Hilbert norm on vector fields. For this purpose, let $V$ denote a Hilbert of vector fields on $\mR^d$ that will assumed, in order to recover elements of $\Diff^p_0$ after solving \cref{eq:ode}, to be continuously included in the space $C^p_0(\mR^d, \mR^d)$ of $C^p$ vector fields that vanish at infinity. This means that $V\subset C^p_0(\mR^d, \mR^d)$ and that, for some constant $c$, one has (letting $\|\cdot\|_\infty$ denote the supremum norm)
\[
\sum_{k=0}^p \|d^k v\|_\infty \leq c \|v\|_V.
\]
To satisfy this assumption, $V$ can be built as a Hilbert Sobolev space of high enough order. In addition, since the continuous inclusion implies that $V$ is a reproducing kernel Hilbert space (RKHS) of vector fields, RKHS theory can be used to build a large variety of Hilbert spaces of interest that satisfy the inclusion property \citep{aronszajn_theory_1950,kadri2016operator,younes2019shapes}. One can then define the action functional of a diffeomorphic motion $((t,x)\in [0,1]\times \mR^d \mapsto \phi(t,x)\in \mR^d)$ as
\[
\int_0^1 \|v(t, \cdot)\|^2_V dt
\]
with $\prt_t \phi(t,x) = v(t, \phi(t,x))$.
A geodesic diffeomorphic motion is an extremal of this action functional and a minimizing geodesic motion minimizes the functional subject to fixed boundary conditions at $t=0$ and $t=1$. In particular, the geodesic distance between two diffeomorphisms $\phi_0$ and $\phi_1$ is defined as 
\begin{multline}
\label{eq:dist.diff}
    d_V(\phi_0, \phi_1) =\\
\inf\left\{ \left(\int_0^1 \|v(t, \cdot)\|^2_V dt\right)^{1/2}:\prt_t \phi(t,x) = v(t, \phi(t,x)), \phi(0, \cdot) = \phi_0, \phi(1, \cdot) = \phi_1\right\} .
\end{multline}

Note that the set over which the infimum is computed may be empty, in which case the distance is infinite. If this set is not empty, then one says that $\phi_1$ is attainable from $\phi_0$. Diffeomorphisms that are attainable from the identity form a subgroup of $\Diff^p_0$, denoted $\Diff_V$, and this subgroup is complete for the geodesic distance \citep{trouve1995action,younes2019shapes}. (Because not every diffeomorphism in $\Diff^p_0$ is attainable from the identity, one is actually building a sub-Riemannian metric on this space. See \citet{arguillere_shape_2014,younes_sub-riemannian_2020}.)

By construction, the distance is right-invariant, i.e., 
\[
d_V(\phi_0, \phi_1) = d_V(\id{d}, \phi_1\circ \phi_0^{-1}),
\]
and this implies that it can be used to define a distance on $\CS$ via
\begin{align*}
d_{\CS}([m_0], [m_1]) &= 
\inf\left\{d_V(\id{d}, \phi): [\phi\cdot m_0] = [m_1] \right\} \\
& = \inf\left\{d_V(\id{d}, \phi): \phi\cdot m_0 \in [m_1]\right\}.
\end{align*}
The distance on $\CS$ can itself be defined directly as
\begin{multline}
\label{eq:shape.oc}
    d_V([m_0], [\phi_1]) =\\
\inf\left\{ \left(\int_0^1 \|v(t, \cdot)\|^2_V dt\right)^{1/2}:\prt_t m(t,\cdot) = v(t, m(t,\cdot), m(0, \cdot) = m_0, m(1, \cdot) \in [m_1]\right\} .
\end{multline}
This provides an optimal control problem in $\CS$ where the control is the time-dependent vector field $v$, and the state equation the ODE $\prt_t m(t,\cdot) = v(t, m(t,\cdot))$. The optimal trajectory transforms the initial $m_0$ into an embedding that is a reparametrization of $m_1$ and provides a minimizing geodesic in $\CS$. If the Sobolev inclusion discussed above holds for $p\geq 1$ at least, the variational problems described in \cref{eq:dist.diff} and \cref{eq:shape.oc} are well defined. The condition that  $\int_0^1 \|v(t, \cdot)\|^2_V dt<\infty$ implies that solutions to the state equations ($\prt_t \phi = v\circ \phi$ or $\prt_t m = v\circ m$) exist and are unique (given initial conditions) over the full unit time interval, ensuring that the optimal control problem is well specified. Moreover, as soon as $\phi_1$ (resp. $[m_1]$) is attainable from $\phi_0$ (resp. $[m_0]$), an optimal solution to the considered problem always exists. Finally, under very mild assumptions on initial conditions, solutions of the geodesic equations exist and are uniquely specified by their initial position and velocity, i.e., $m(0, \cdot)$ and $\prt_t m(0, \cdot)$ for $d_\CS$. The geodesic equation is the Euler-Lagrange equation associated with the variational problem, satisfied by stationary points of \cref{eq:dist.diff} or \cref{eq:shape.oc} (equivalently, they are the equations provided by Pontryagin's maximum principle). In the case considered here, they are special instances of the geodesic equations for right-invariant Riemannian metrics on Lie groups, as described in \citet{arnold1966sur,arnold1978mathematical} and are often referred to as Euler-Arnold equations \citep{arnold2021topological} or Euler-Poincar\'e equations \citep{ebin1970groups,holm1998euler--poincare}.

In practice, one does not solve this problem exactly, but relaxes the endpoint condition $m(1,\cdot) \in [m_1]$ by adding a penalty term, therefore minimizing 
\begin{equation}
    \label{eq:.lddmm}
\int_0^1 \|v(t, \cdot)\|^2_V dt + U([m(1, \cdot)], [m_1])
\end{equation}
subject to $\prt_t m(t,\cdot) = v(t, m(t,\cdot))$. In many of the applications, the function $U$ takes the form
\[
U([m_0], [m_1]) = \|\mathcal J_{[m_0]} - \mathcal J_{[m_1]}\|_H^2
\]
where $[m] \mapsto \mathcal J_{[m]}$ is a mapping from $\CS$ into a (much larger) Hilbert space $H$. These ``chordal metrics'' use representations of embedded curves or surfaces as measures, currents or varifold. For simplicity the presentation below will ignore this relaxation step (which is however necessary to make the computation numerically feasible) and work as if the endpoint conditions are exact. The reader is referred to \citet{bauer2019metric} or \citet{charon2020fidelity}, and to the references within, for more information on chordal metrics. 

\begin{figure}
    \centering
    \includegraphics[trim=3cm 1cm 3cm 1cm,clip,width=0.2\textwidth]{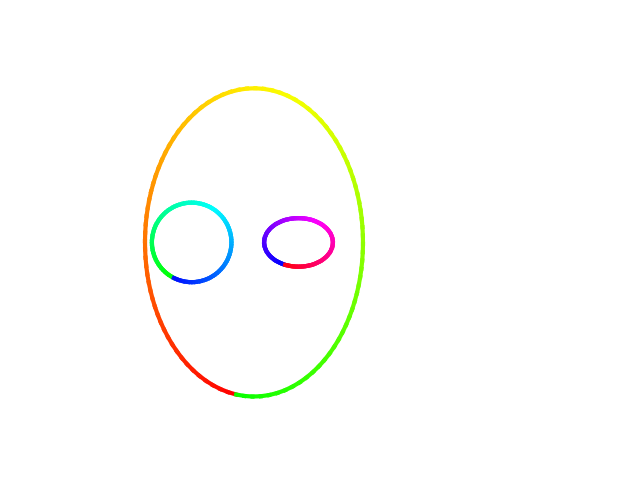}
    \includegraphics[trim=3cm 1cm 3cm 1cm,clip,width=0.2\textwidth]{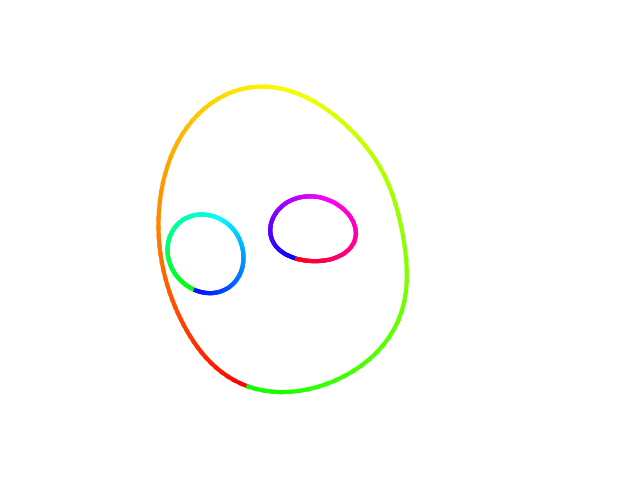}
    \includegraphics[trim=3cm 1cm 3cm 1cm,clip,width=0.2\textwidth]{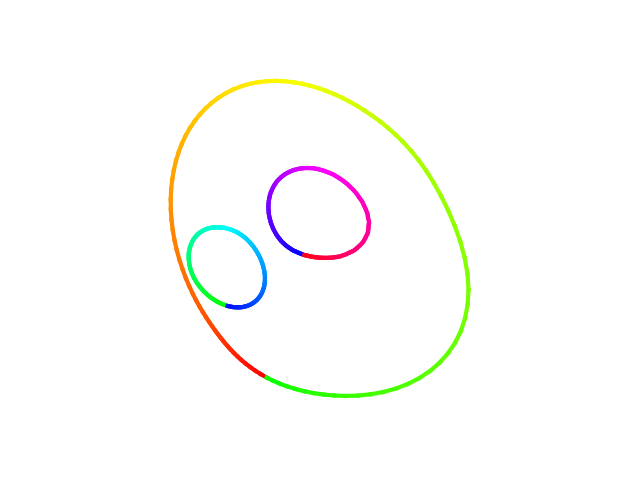}
    \includegraphics[trim=3cm 1cm 3cm 1cm,clip,width=0.2\textwidth]{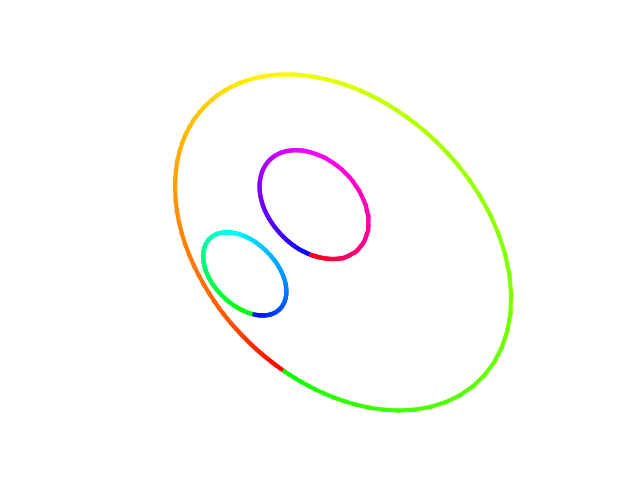}
    \caption{Four time points of a geodesic evolution in shape space. Note that shapes in this example have multiple components. Contour coloring match across time points and track the evolution of the curve initial parametrization.}
    \label{fig:geod.curves}
\end{figure}

A two-dimensional example of geodesic is presented in \cref{fig:geod.curves}. These geodesics provide the non-linear equivalent of a linear interpolation in Euclidean space. 

\section{Hybrid models}
\label{sec:hybrid}

\subsection{Description}
The previous framework can be slightly extended to allow the norm used in the shape space to depend on the shape itself, replacing the control cost in \cref{eq:shape.oc} by
\[
\int_0^1 \|v(t, \cdot)\|^2_{[m(t, \cdot)]} dt,
\]
so that the cost depends on both control and state. This still provides a sub-Riemannian distance in shape space, and the problem remains well specified as soon as one ensures that the shape-dependent norms still control the norm on $V$, so that an inequality ensuring
\[
\|v\|_V \leq C \|v\|_{[m]}
\]
holds for all $m\in \CS$ and $v\in V$ (where the upper-bound may be infinite). Typical applications of this construction use a ``weak norm'' $v \mapsto \llbracket v \rrbracket_{[m]}$ (which, by itself would not guarantee the existence of solutions to the state equation), possibly motivated by material or biological constraints, ``regularized'' by the norm on $V$, therefore taking
\begin{equation}
    \label{eq:hybrid}
\|v\|_{[m]}^2 = \kappa \|v\|_V^2 + \llbracket v \rrbracket_{[m]}^2
\end{equation}
for some $\kappa>0$.

The following section discusses several possible choice for $\llbracket v \rrbracket_{[m]}$ in \cref{eq:hybrid}, in which the shape is considered as an elastic material and the norm corresponds to the elastic energy associated with an infinitesimal displacement along $v$ (the reader may refer to, e.g.,  \citet{ciarlet1988three-dimensional,gonzalez2008first} for more details on elasticity concepts that are used below). The concept of  ``hybrid'' metrics in \cref{eq:hybrid} was suggested in \citet{younes2018hybrid}. A similar approach for spaces of images (combined with a ``metamorphosis'' metric \citep{miller2001group,trouve2005metamorphoses}) was introduced in \citet{berkels2015time}, and metrics formed as discrete iterations of small elastic deformations were also studied in \citet{wirth2011continuum}.

\subsection{Elastic metrics}
\subsubsection{Three-dimensional case}
\label{sec:elastic.3D}
The energy of a hyper-elastic material $\Omega$ subject to a deformation $\varphi$ takes the form (letting $\Id{d}$ denote the identity matrix in $\mR^d$)
\[
E = \int_{\varphi(\Omega)} G(x, \varphi(x))  dx
\]
where
\[
G(x, \varphi) = W\left(x, d\varphi^T d\varphi - \Id{3}\right)
\]
for a function  $W: \Omega \times \Sym^+ \to [0, +\infty)$ (where $\Sym^+$ is the set of 3 by 3 positive semi-definite matrices) such that $W(x,S) = 0$ if and only if $S=0$. The matrix $C = d\varphi^T d\varphi$ is the Cauchy-Green strain tensor, which is such that $u^T C u = |d\varphi\, u|^2$, and $W$ measures the deviation of this tensor from the identity matrix.

A second-order expansion of $G$, near $\varphi = \id{3}$ takes the form (using the fact that $\partial_2W(x, 0) = 0$)
\begin{equation}
    \label{eq:elastic.order2}
G(x, \varphi) \simeq \frac12 \partial_2^2W(x, 0)(dv + dv^T, dv + dv^T)
\end{equation}
where $v = \varphi - \id{3}$. Here $\partial_2 W(x, 0)$ and $\partial_2^2 W(x, 0)$ are the first and second derivative with respect to the second variable of $W$, therefore a positive semi-definite symmetric bilinear form on $\Sym$ (the space of 3 by 3 symmetric matrices). 

This can be used to define an  elastic metric on 3D vector fields. Here, $\Omega$ is considered as an ``unparametrized shape,'' taking the role of $[m]$ in the previous sections. Using the previous notation, this corresponds to taking the manifold $M$ to be an open subset of $\mR^3$ (e.g., an open ball), $m$ an embedding of $M$ into $\mR^3$ and identifying $\Omega = m(M)$ to $[m]$. 
The hybrid norm will therefore be denoted 
\[
\|v\|^2_\Omega = \kappa \|v\|_V^2 + \llbracket v\rrbracket_\Omega^2
\]
and the rest of the discussion focuses on $\llbracket v\rrbracket_\Omega$.
Based on \cref{eq:elastic.order2}, one is led to define a 3D elastic metric on vector fields as any norm taking the form
\[
\llbracket v \rrbracket_{\Omega}^2 = \int_{\Omega} B(x, \epsilon(x)) dx
\]

where $\varepsilon(x) = (dv(x) + dv(x)^T)/2$ is known as the infinitesimal strain tensor of the deformation and $B(x, \cdot)$ is a positive semi-definite quadratic form on $\Sym$, typically referred to as the elastic tensor. Generically, $B(x, \cdot)$ can be  represented as a $6\times 6$ symmetric positive semi-definite matrix, that is with 21 parameters in total at each $x$. In a majority of applications however, model symmetry assumptions significantly reduce the complexity of this elasticity tensor. In particular, in the case of a uniform and isotropic material, $B(x,\cdot)$ is independent of the position and takes the specific form
\begin{equation}
\label{eq:isotropic} 
B(x, \varepsilon) = B(\varepsilon) = \frac{\lambda}{2} \trace(\epsilon)^2 + \mu \trace(\epsilon^2)
\end{equation}

which is the linearization of the energy of a Saint Venant-Kirchhoff material. In that case, the elasticity tensor is only described by the two parameters $\lambda$ and $\mu$ which are called the Lam\'e coefficients of the material. 

To provide another example, consider the case of a partially isotropic and laminar model, introduced in \citet{hsieh2019model,hsieh2021diffeomorphic,hsieh2022mechanistic} under the assumption that $\Omega$ can be parametrized by a foliation.  More precisely, assume that there exist two surfaces  ${\mathcal M}_{\mathrm{bottom}}$ and ${\mathcal M}_{\mathrm{top}}$ (bottom and top layers) included in $\partial\Omega$ and  a diffeomorphism $\varPhi: [0, 1] \times {\mathcal M_{\mathrm{bottom}}} \rightarrow \Omega$ such that $\varPhi(\{0\}\times {\mathcal M_{\mathrm{bottom}}}) = {\mathcal M_{\mathrm{bottom}}}$ and $\varPhi(\{1\} \times {\mathcal M_{\mathrm{bottom}}}) = {\mathcal M_{\mathrm{top}}}$. Let $\mathcal M_s = \varPhi(\{s\}\times {\mathcal M_{\mathrm{bottom}}})$, $s \in [0, 1]$, denote ``the layer at level $s$,'' $S$ the  transverse vector field $S = \partial_s \varPhi$ and $N$ a unit vector field normal to all $\mathcal M_s$. One then introduces the following elasticity tensor:
\begin{align}
	\label{eq:computation_elastic_energy}
	B(x, \varepsilon)=&\lambda_{\mathrm{tan}} \left( \trace(\varepsilon) - N^T \varepsilon N \right)^2
		+\mu_{\mathrm{tan}}\left(\trace(\varepsilon^2) - 2 \, N^T \varepsilon^2 N + (N^T \varepsilon N)^2
		\right)
		\\ \nonumber
		&+ \mu_{\mathrm{tsv}} \, (S^T \varepsilon S)^2
			+ 2 \, \mu_{\mathrm{ang}} \left( S^T \varepsilon^2 S - (N^T \varepsilon S)^2 \right) ,
\end{align}
The first two terms in this expression define an isotropic model on each layer. The third term measures a transversal string, evaluated along $S$. The last term measures an angular strain, that vanishes when $S$ is normal to the layers. Here, the coefficients $\lambda_{\mathrm{tan}}, \mu_{\mathrm{tan}}, \mu_{\mathrm{tsv}}, \mu_{\mathrm{ang}}$ must be constant on each layer $m_s$ (they may depend on $s$). Note that, if $\tau_1, \tau_2$ are two orthonormal vectors fields that are tangent to the layers so that $(\tau_1, \tau_2, N)$ forms at all points an orthonormal frame, then 
\[
\trace(\varepsilon) - N^T \varepsilon N = \tau_1^T \varepsilon \tau_1 + \tau_2^T \varepsilon \tau_2
\]
and 
\[
\trace(\varepsilon^2) - 2 \, N^T \varepsilon^2 N + (N^T \varepsilon N)^2 = (\tau_1^T \varepsilon \tau_1)^2 + (\tau_2^T \varepsilon \tau_2)^2 + 2(\tau_1^T \varepsilon \tau_2)^2
\]
so that the first two terms only involve deformations tangent to the layers.

Importantly, the space of such ``layered structures'' is stable by diffeomorphic action. Indeed, given $\Omega$ and $\varPhi$ as above, and $\varphi$ a diffeomorphism of $\mR^3$, one defines the transformed structure by:
\[
\varphi\cdot (\Omega, \varPhi) = (\varphi(\Omega), \varphi \circ \varPhi \circ \varphi^{-1}).
\]
In particular, $S$ transforms through $\varphi$ as $\phi\cdot S = (d\phi S) \circ \phi^{-1}$.

\medskip

Returning to the general case, one must emphasize the fact that the action functional
\[
\int_0^1 \int_{\Omega(t)} B(x, \varepsilon(x)) dx dt
\]
with $\prt_t \varphi(t,x) = v(t, \varphi(t,x))$ and $\Omega(t) = \varphi(t,\cdot)(\Omega_0)$ is not the energy of a deforming elastic material, in the sense given to it in elasticity theory. In contrast, it may be understood as a sum of infinitesimal elastic energies, for a volume that slowly deforms, and at each time step, remodels its structure to reach an equilibrium state {\em without---up to reorientation---changing its elasticity properties}.

\subsubsection{Elastic metrics on surfaces}
\label{sec:surfaces}
The definition of elastic metrics on surfaces can be inferred using a pattern similar to the 3D derivation. Let $\mathcal M$ be a surface in $\mR^3$ and consider a one-to-one immersion $\phi: \mathcal M \to \mR^3$ (one can, in this discussion, think of $\phi$ as the restriction to $\mathcal M$ of a diffeomorphism of $\mathbb R^3$). To define a hyperelastic energy, assume (restricting $\mathcal M$ if needed and introducing partitions of unity) that two vector fields $\tau_1, \tau_2$ are chosen on $\mathcal M$ such that they form at each point an orthonormal frame, and let $\nu = \tau_1\times \tau_2$. Let $F(x)$ denote the $3\times 2$ matrix $[d\phi\, \tau_1, d\phi\, \tau_2](x)$, where the 3D columns are expressed in the canonical basis of $\mR^3$, and consider energies of the form
\[
\int_M W(x, F(x)) d\mathrm{vol}_M (x).
\]
Material independence requires that $W$ is invariant when $F$ is multiplied on the right by a 2D rotation matrix, and this implies that $W$ only depends on $FF^T$. To obtain the expression of the metric, we let $\phi(x) = x + v(x)$ and make a first order expansion in $v$ of $FF^T$ with $F = [\tau_1 + dv \tau_1, \tau_2 + dv \tau_2]$ yielding
\[
FF^T \simeq  \pi_{\mathcal M} + \pi_{\mathcal M} dv^T + dv \pi_{\mathcal M}
\]
where $\pi_{\mathcal M} = \tau_1\tau_1^T + \tau_2\tau_2^T$ is the orthogonal projection on the tangent plane to $\mathcal M$ at $x$. The Riemannian elastic metric should therefore be taken as a quadratic form of $\eta_{\mathcal M} := (\pi_{\mathcal M} dv^T + dv \pi_{\mathcal M})/2$. Expressing this operator in the basis $(\tau_1, \tau_2, \nu)$, one sees that it depends on the five quantities $a_{11} = \tau_1^T dv \tau_1$, $a_{22} = \tau_2^T dv \tau_2$, $a_{12} + a_{21} = \tau_1^T dv \tau_2 + \tau_2^T dv \tau_1$, $a_{13} = \nu^T dv \tau_1$ and $a_{23} = \nu^T dv \tau_2$, yielding 15 free parameters for the ``elastic norm'' $\llbracket v \rrbracket_{\mathcal M}$. (Like in the previous section, an identification is made between the unparametrized surface $\mathcal M = m(M)$ and the equivalence class $[m]$.)

The norm is isotropic if it satisfies $\llbracket Rv \rrbracket_{\mathcal M} = \llbracket v \rrbracket_{\mathcal M}$ for any 3D rotation that leaves $\nu$ invariant. This implies that the matrix $\mathbf a = \begin{pmatrix} a_{11} & (a_{12} + a_{21})/2 \\ (a_{12} + a_{21})/2 & a_{22} \end{pmatrix}$ is transformed by a 2D rotation as $\mathbf a \mapsto R^T \mathbf a R$ and the vector $\mathbf b = \begin{pmatrix} \nu^T dv \tau_1 \\ \nu^T dv \tau_2\end{pmatrix}$ as $\mathbf b\mapsto \mathbf b R$. Using usual invariance arguments, this requires that  the squared norm must be a (quadratic) function of $\mathrm{trace}(\mathbf a)$, $\mathrm{trace}(\mathbf a^2)$ and $|\mathbf b|^2$, yielding
\begin{equation}
    \label{eq:surf.elastic}
\llbracket v \rrbracket_{\mathcal M}^2 = \int_{\mathcal M} \beta(x, \eta_{\mathcal M}) d\mathrm{vol}(x)
\end{equation}
with
\begin{equation}
\label{eq:surf.elastic.metric}
\beta(x, \eta_{\mathcal M}) = \lambda_{\mathrm tan} \mathrm{trace}(\mathbf a)^2 + \mu_{\mathrm tan} \mathrm{trace}(\mathbf a^2) + \mu_{\mathrm{tsv}} |\mathbf b|^2.
\end{equation}
The three different terms of this metric can be also interpreted as penalties on the changes of local area, metric tensor and normal vector respectively, as pointed out in \citet{jermyn2012elastic} (see also their intrinsic expressions derived in \cref{app:laminar}). \citet{jermyn2012elastic} focuses on the special case $\lambda_{\mathrm tan} = 1/16, \mu_{\mathrm tan} = 0, \mu_{\mathrm{tsv}}=1$, which can be shown to be isometric to a Euclidean metric under a ``square root normal transform.''

To consider another example, let $\lambda_{\mathrm tan} = 0$ and $\mu_{\mathrm tan} = \mu_{\mathrm{tsv}} = 1$. Then
\begin{align}
\nonumber
\beta(x, \eta_{\mathcal M}) &= a_{11}^2 + a_{22}^2 + \frac12 (a_{12} + a_{21})^2  + a_{13}^2 + a_{23}^2\\
\nonumber
&= a_{11}^2 + a_{22}^2 + a_{12}^2 + a_{21}^2   + a_{13}^2 + a_{23}^2 - \frac12 (a_{12} - a_{21})^2\\
\label{eq:h1.corrected}
&= \trace(dv dv^T) - \frac12 (a_{12} - a_{21})^2
\end{align}
The first term, $\trace(dv dv^T)$, corresponds to the $H^1$ metric on $\mathcal M$, used e.g., in \citet{younes2018hybrid}. This metric, without the correction term $\frac12 (a_{12} - a_{21})^2$ is not an elastic metric. It belongs however to a larger class of metrics, studied in \citet{su2020shape}, where the correction term is added to \cref{eq:surf.elastic.metric} with a fourth parameter.

\begin{figure}
    \centering
    \includegraphics[trim=2cm 2cm 10cm 2cm, clip, width=0.3\textwidth]{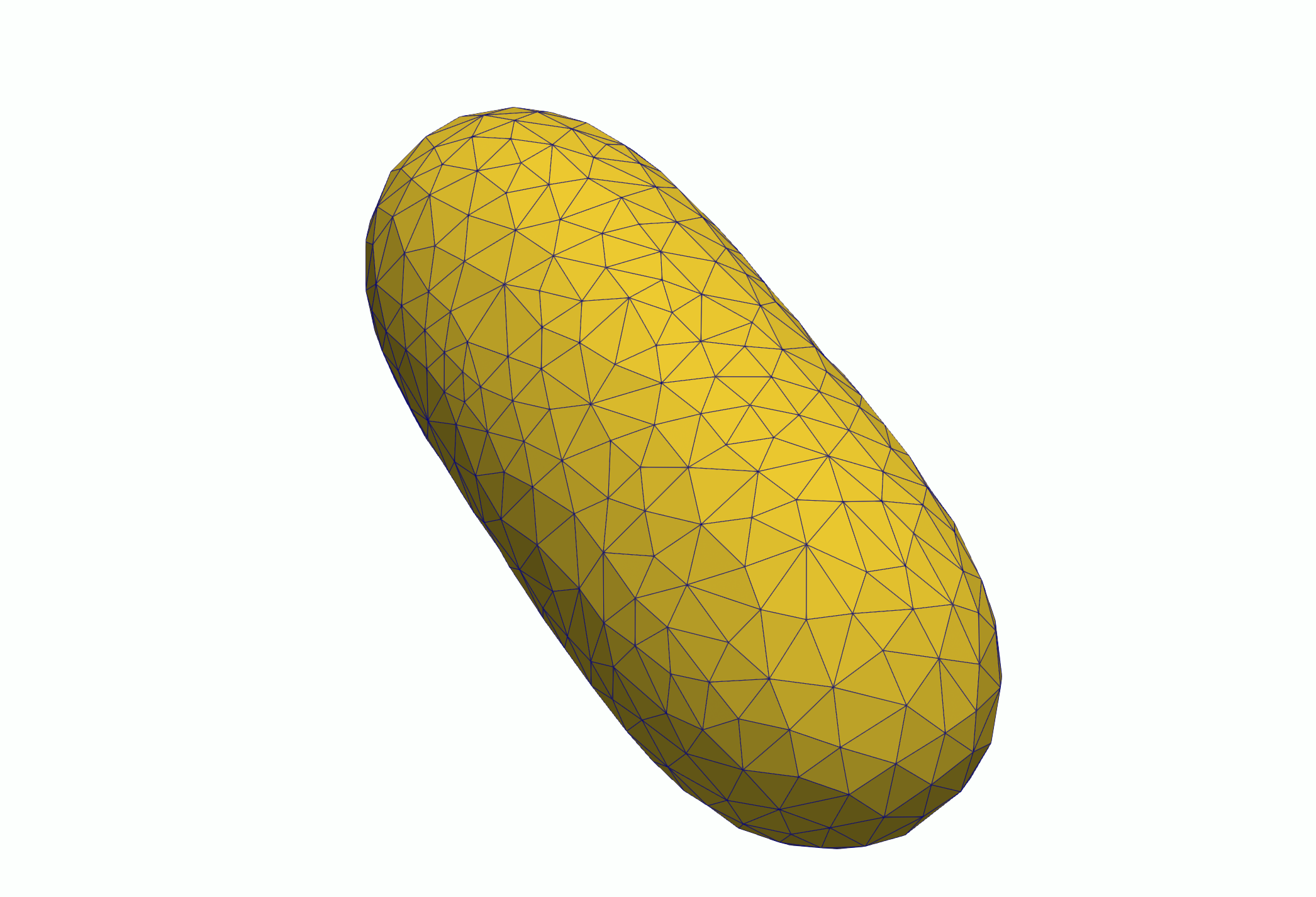}
    \hskip 0.2\textwidth
    \includegraphics[trim=2cm 2cm 10cm 2cm, clip,width=0.3\textwidth]{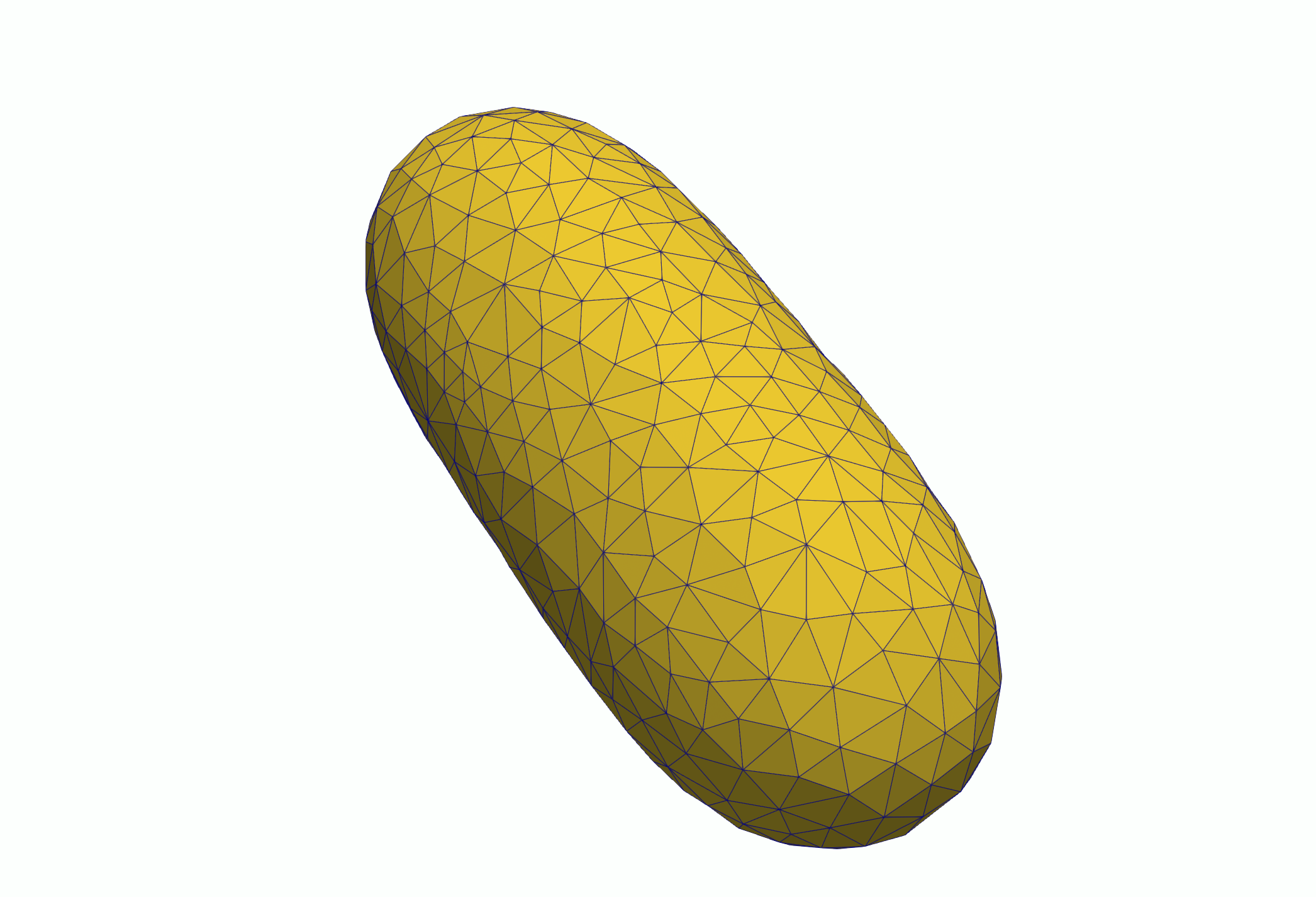}\\
    \includegraphics[trim=2cm 2cm 10cm 2cm, clip,width=0.3\textwidth]{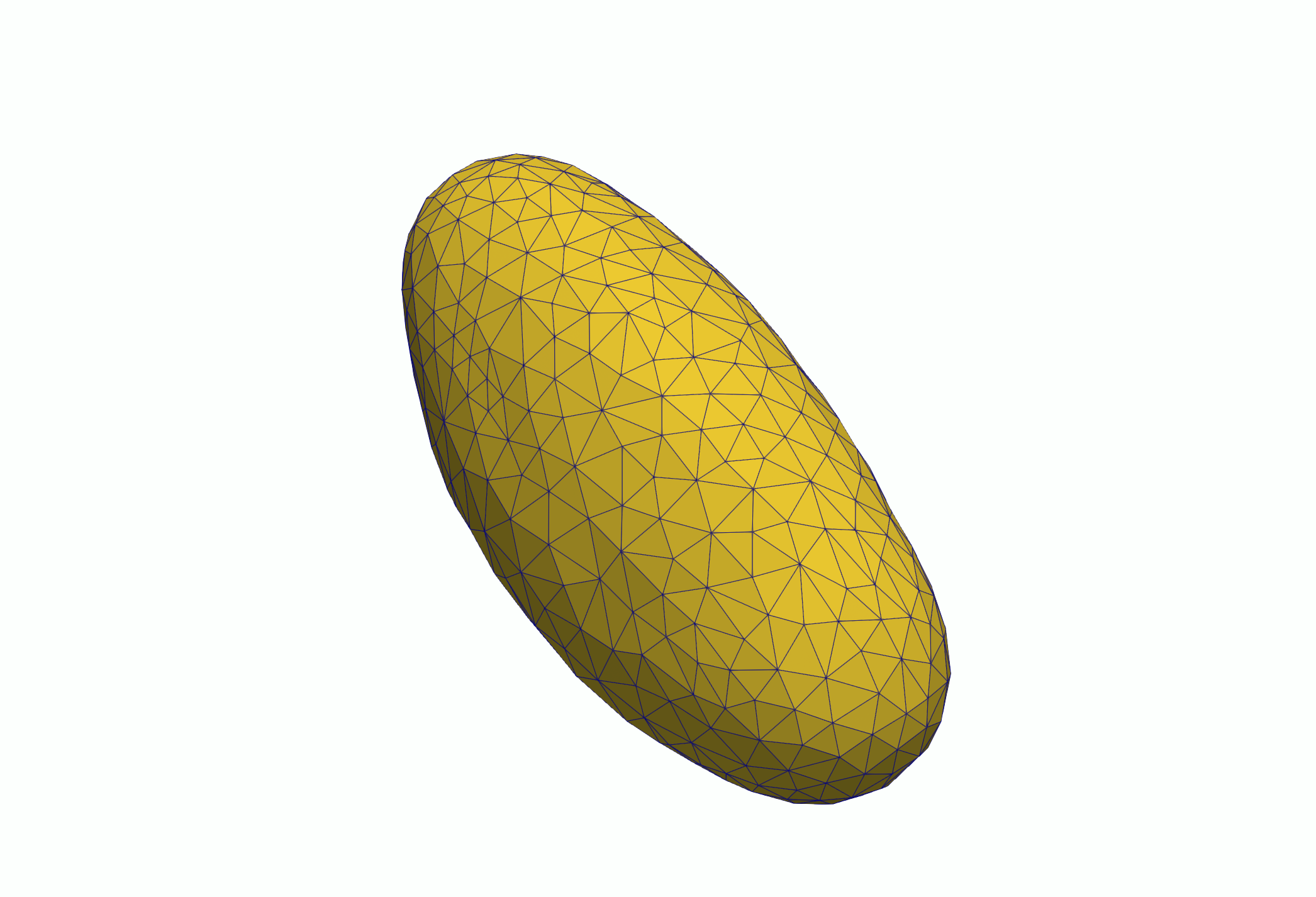}
    \hskip 0.2\textwidth
    \includegraphics[trim=2cm 2cm 10cm 2cm, clip,width=0.3\textwidth]{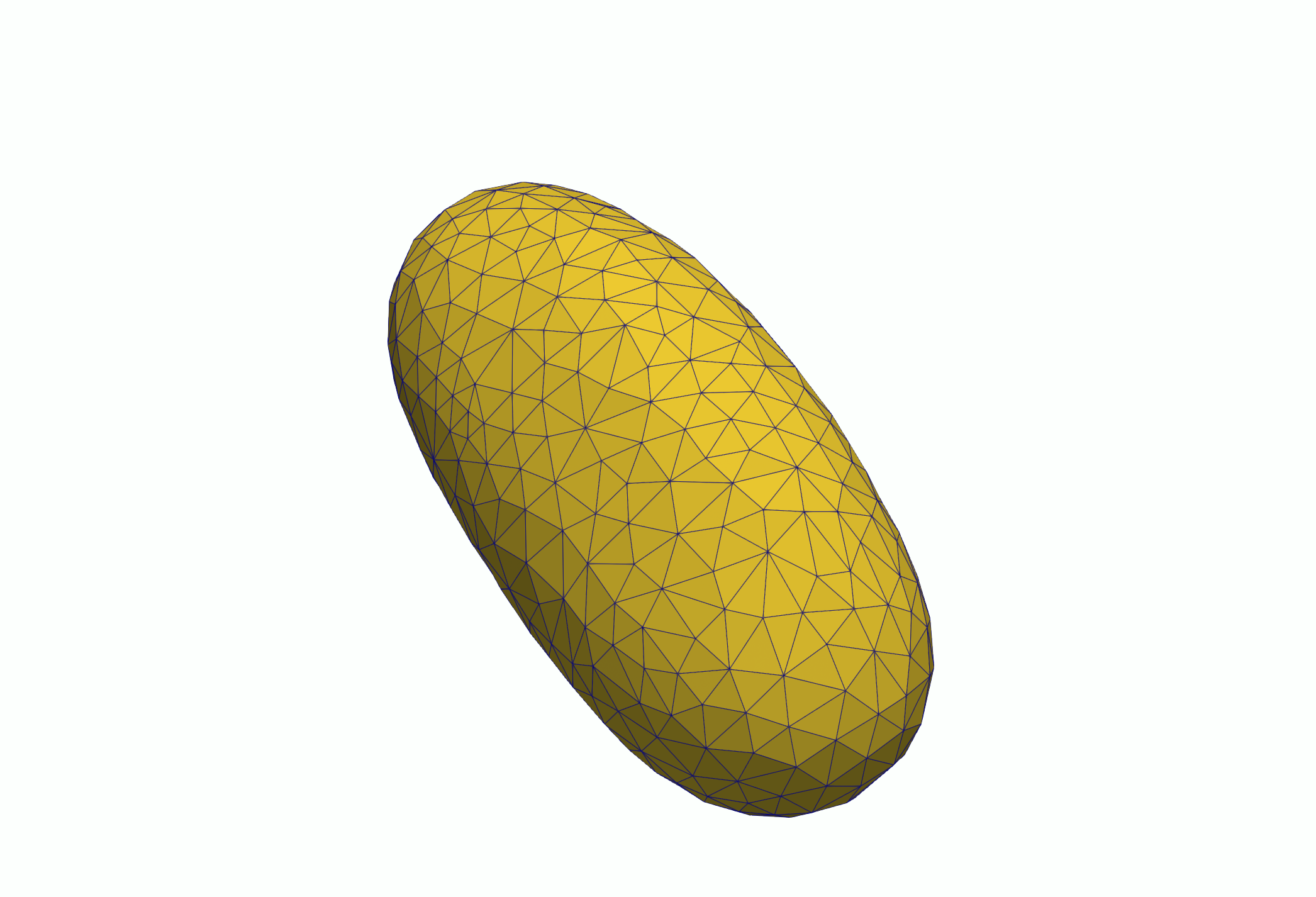}\\
    \includegraphics[trim=2cm 2cm 10cm 2cm, clip,width=0.3\textwidth]{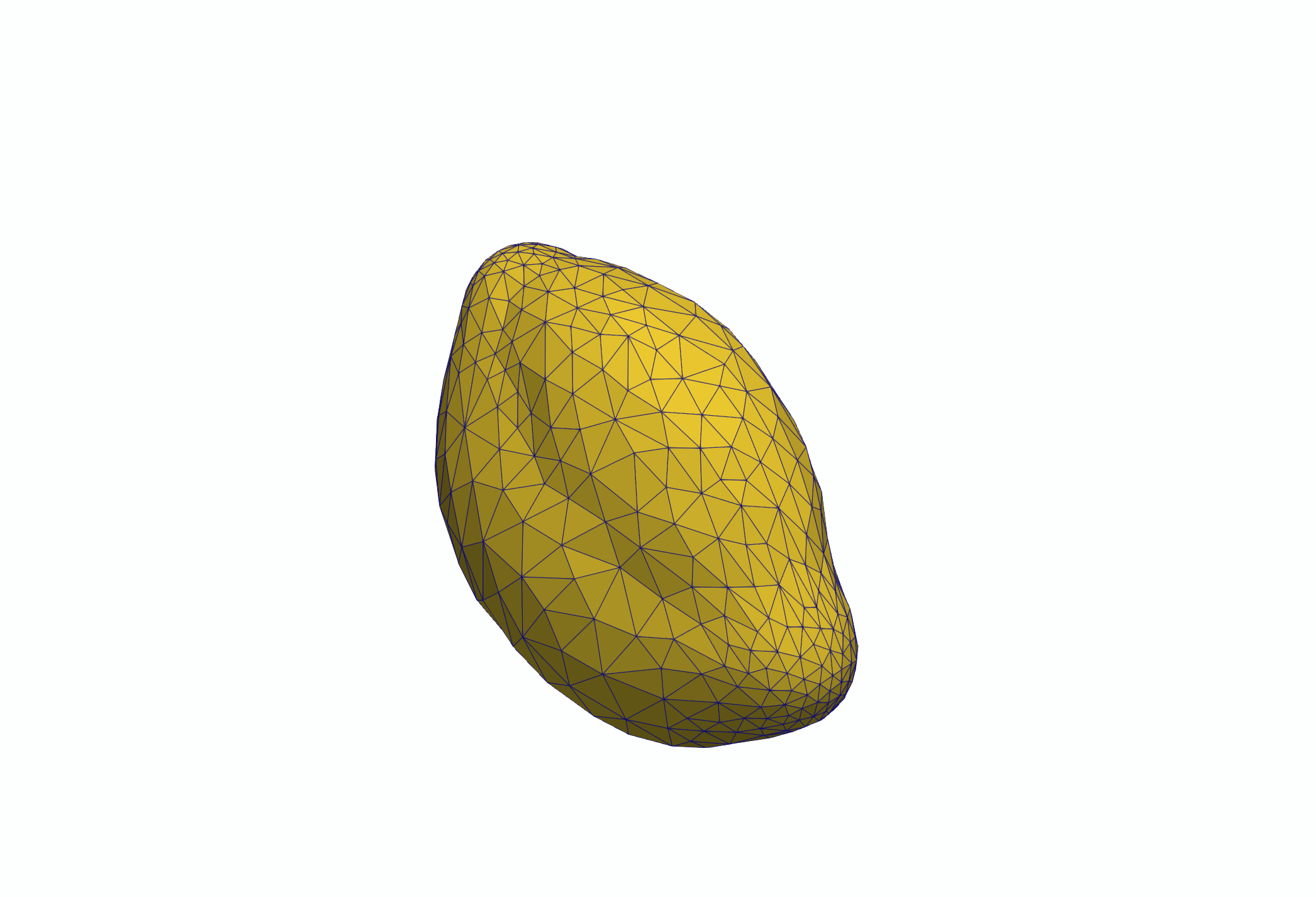}
    \hskip 0.2\textwidth
    \includegraphics[trim=2cm 2cm 10cm 2cm, clip,width=0.3\textwidth]{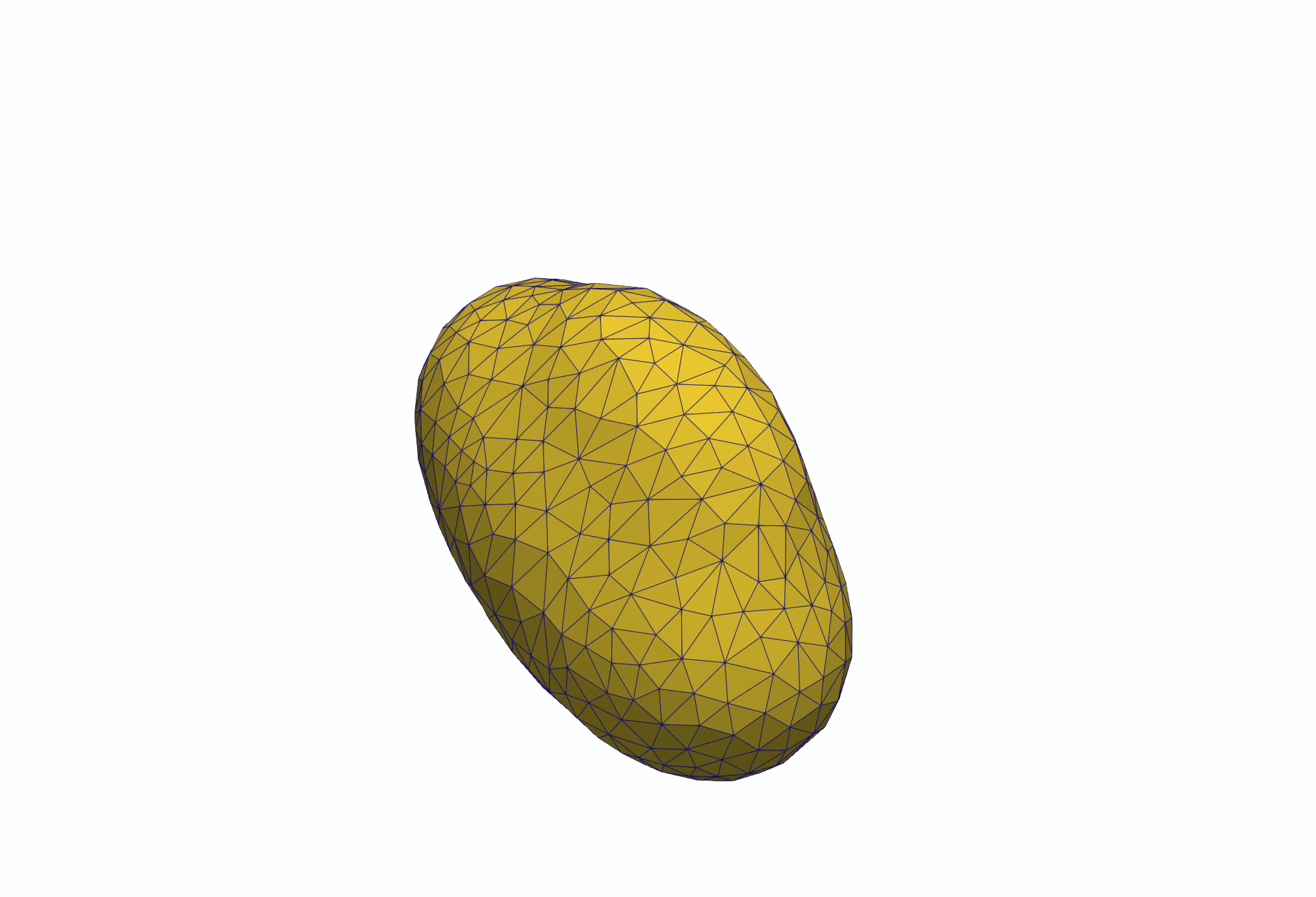}\\
    \includegraphics[trim=2cm 2cm 10cm 2cm, clip,width=0.3\textwidth]{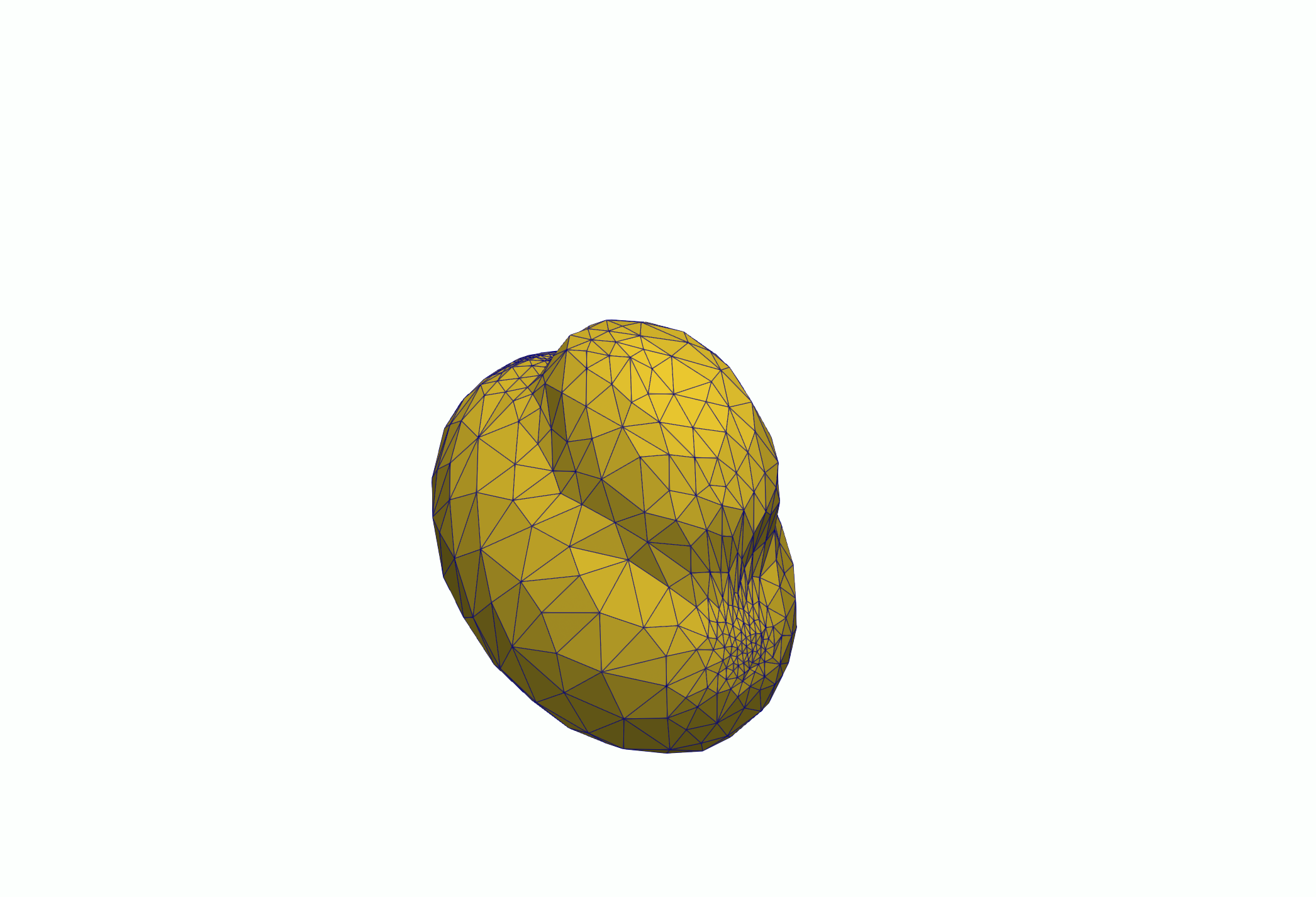}
    \hskip 0.2\textwidth
    \includegraphics[trim=2cm 2cm 10cm 2cm, clip,width=0.3\textwidth]{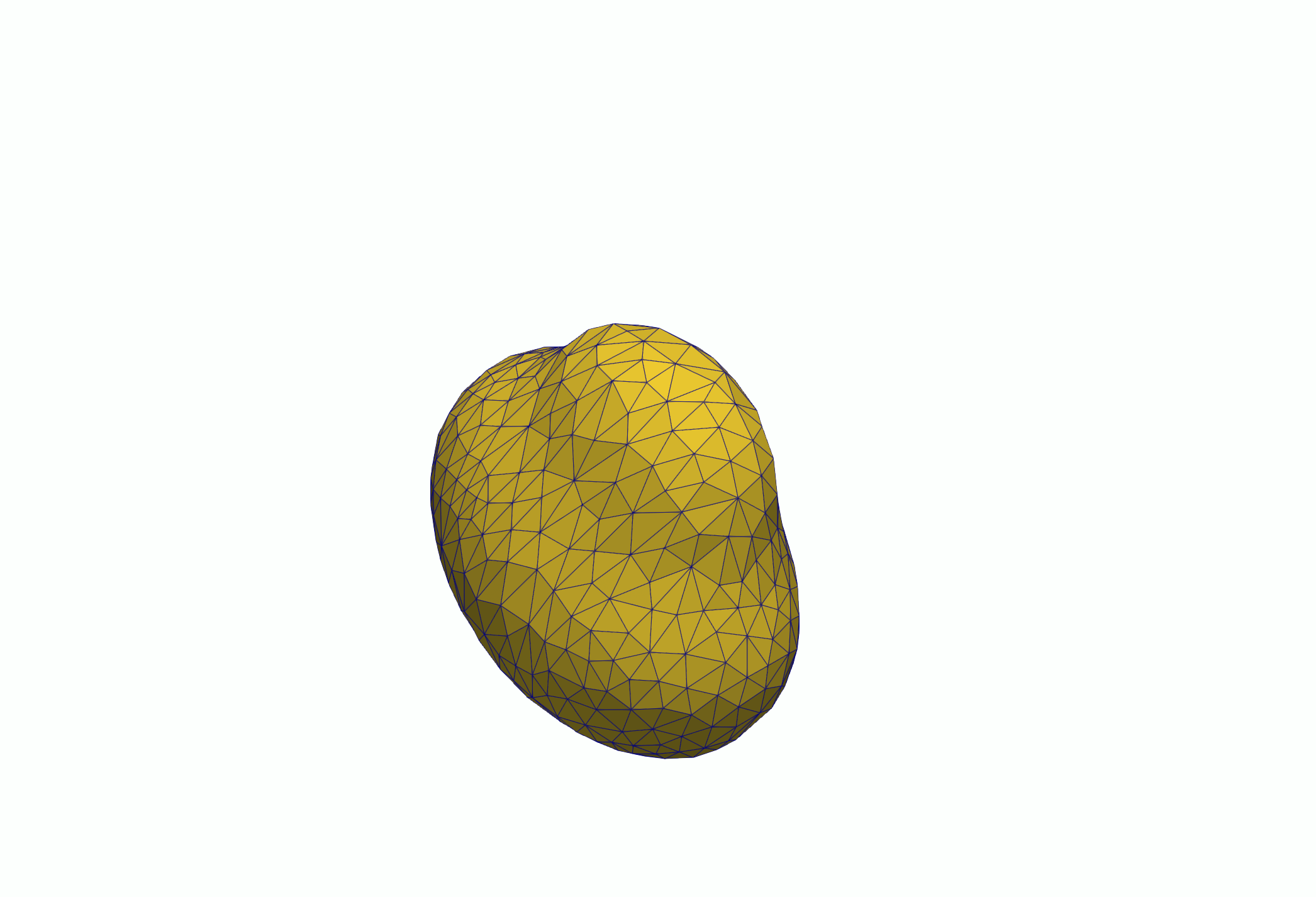}
    \caption{Comparison between geodesics between surfaces using a pure LDDMM and a hybrid LDDMM/elastic metric. First column: Four time points of an LDDMM geodesic ($t=0$, $t=0.3$, $t=0.7$ and $t=1$). Second column: same time points for the hybrid geodesic. One can note a difference in the intermediate shapes, and (as indicated by the triangulation) higher local contraction associated with the LDDMM metric. The hybrid metric uses the expression provided in \cref{eq:h1.corrected}.}
    \label{fig:elastic.surfaces}
\end{figure}

Such elastic metrics can be used in combination with the LDDMM metric through the hybrid setup described above. As an illustration, \cref{fig:elastic.surfaces} provides a comparison of the geodesic trajectories between two surfaces, obtained with the pure LDDMM model and a hybrid model using the elastic term given by \cref{eq:h1.corrected}.

\medskip

The norm in \cref{eq:surf.elastic} can also be obtained as a limit of the laminar elastic model of the previous section and the energy in \cref{eq:computation_elastic_energy}, which is shown in \cref{app:laminar} by also providing an intrinsic expression of the elastic norm. This in part justifies the terminology of elastic metrics given to this framework in the related literature.

% we show that  leads, in the limit of infinitely thin shell materials, to a class of first-order invariant metrics on the space of surfaces which have played an important role in shape analysis since their introduction in works such as \cite{younes1998computable,srivastava2016functional,su2020shape}. This justifies in particular the appellation of elastic metrics given to those in this literature.

\subsubsection{Elastic metrics on curves}
If $\mathcal M$ is a 3D curve, the same analysis shows that elastic metrics should depend on the products $\tau^T dv \tau$, $\nu_1^T dv \tau$ and $\nu_2^T dv \tau$, where $\tau$ is a unit tangent on $M$ and $(\tau, \nu_1, \nu_2)$ is a continuous positively oriented frame defined along the curve. Denote $\prt_s v = dv\tau$ for the derivative with respect to arc length, as introduced, e.g., in \citet{michor2007overview}.  The metric must also be invariant to rotations of the normal frame $(\nu_1, \nu_2)$ and changes of orientation on $\mathcal M$, which requires the metric to take the form
\begin{equation}
    \label{eq:curve.elastic}
\llbracket v \rrbracket_{\mathcal M}^2 = \int_{\mathcal M} \beta(x, \prt_s v) d\mathrm{vol}(x)
\end{equation}
with
\[
\beta(x, \prt_s v)= \mu_{\mathrm tan} (\tau^T \prt_s v)^2 + \mu_{\mathrm{tsv}} ((\nu_1^T \prt_sv)^2 + (\nu_2^T\prt_s v)^2).
\]

\begin{figure}
    \centering
    \includegraphics[trim=3cm 1cm 3cm 1cm,clip,width=0.2\textwidth]{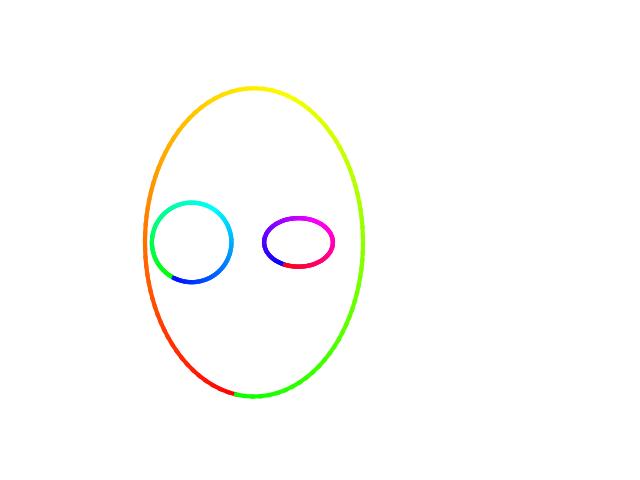}
    \includegraphics[trim=3cm 1cm 3cm 1cm,clip,width=0.2\textwidth]{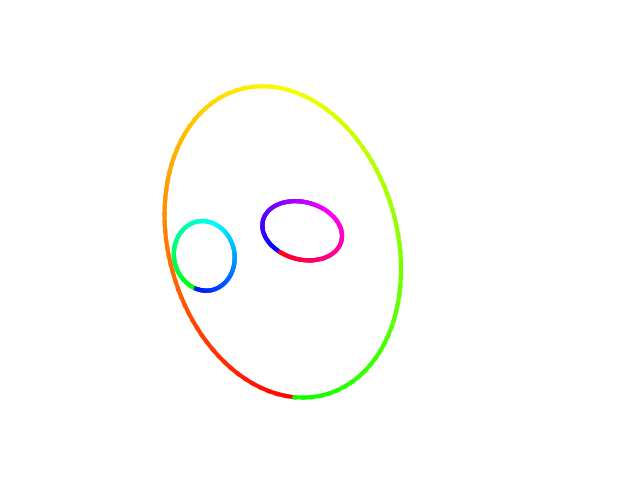}
    \includegraphics[trim=3cm 1cm 3cm 1cm,clip,width=0.2\textwidth]{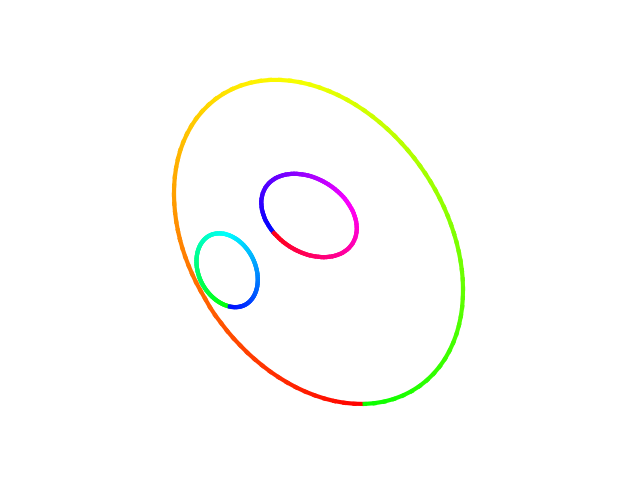}
    \includegraphics[trim=3cm 1cm 3cm 1cm,clip,width=0.2\textwidth]{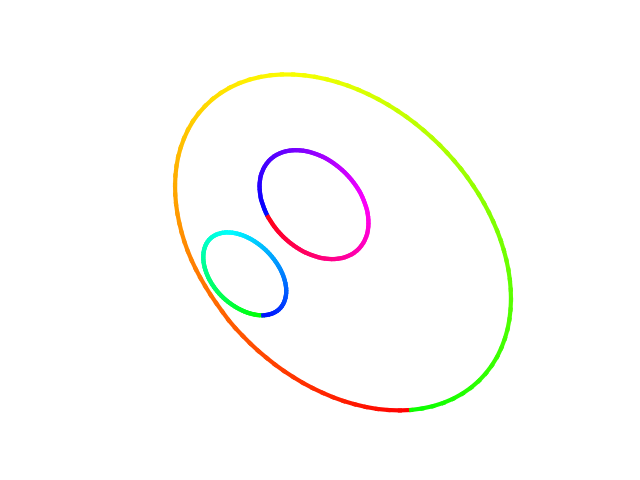}
    \caption{Four time points of a geodesic evolution in shape space for a hybrid metric. The initial and final shapes are the same as those in \cref{fig:geod.curves}, but one can note, in particular, that the elliptical shapes are better conserved during the motion.}
    \label{fig:geod.curves.elastic}
\end{figure}

The special case of planar curves has been extensively discussed. In this case, letting $\nu$ denote the unit normal, the metric has two parameters, with
\[
\beta(x, \prt_s v)= \mu_{\mathrm tan} (\tau^T \prt_s v)^2 + \mu_{\mathrm{tsv}} (\nu^T \prt_s)^2.
\]
When $\mu_{\mathrm tan} = \mu_{\mathrm{tsv}} = 1$, one gets $\beta(x, \prt_s v)= |\prt_s v|^2$. The resulting metric was introduced in \citet{younes1996distance,younes1998computable} and shown to be isometric to a flat metric using a square root transform. This metric was called ``$H^1_0$'' in \citet{mumford2006riemannian} and further studied in \citet{younes2008metric}.  The case $\mu_{\mathrm tan} = 1,  \mu_{\mathrm{tsv}} = 1/4$ was considered in \citet{mio2007shape,srivastava2016functional}, and a similar square root transform was seen to provide an isometry with a flat space in this case also. This isometry was extended to the general case in \citet{younes2018elastic, younes2019shapes} and in \citet{needham2020simplifying} (another isometry was also introduced in \citet{bauer2014constructing}). The reader is referred to the cited references for more details on the exact expression of the isometry. \Cref{fig:geod.curves.elastic} provdes an example of geodesic evolution for a hybrid metric, to be compared with \cref{fig:geod.curves}.

\section{Growth models}
\label{sec:growth}
\subsection{Introduction}
The previous section described various metrics in shape space that are built as a regularized linearized elastic energy. Optimal paths (i.e., geodesics) associated with these metrics  prefer different trajectories from those associated with the ``standard'' spaces discussed in \cref{sec:diffeo}, and  tend to inherit some of the properties suggested by the elastic intuition. However, not all trajectories of interest need to be geodesics for some metric or satisfy a least-action principle.  In particular, including external actions, with in particular possible mechanisms describing growth\footnote{Following common terminology, we consider growth as a general shape change mechanism, also including atrophy, as a ``negative growth.''}, will provide shape analysis methods with additional capability of modeling transformations typically observed in biology or medicine. 

A leading model for shape change in the framework of elasticity theory introduces the notion of {\em morpho-elasticity} in which shapes are subject to the action of a ``growth tensor,'' which partly accounts for the derivative of the deformation, $d\varphi$ (see \citet{goriely2017mathematics} for an extensive introduction to the subject and for references). Letting $G$ denote the growth tensor, one writes $d\varphi = AG$, where $A$ completes the growth tensor to provide a valid differential $d\varphi$, in a way that would minimize the elastic energy (so that one applies the elastic cost to $A^TA - \Id{d}$ rather than to $d\varphi^Td\varphi - \Id{d}$). This approach does not necessarily lead to the trivial solution $A = \Id{d}$ because the growth tensor $G$ is not necessarily ``compatible'', i.e., there may not always exist a transformation $\varphi$ such that
$d\varphi^T d\varphi = G^TG$. 

Considering small deformations, i.e., linearizing $d\varphi = AG$ for $\varphi$ and $G$ close to the identity, and writing $\varphi = \id{d} + v$, $A = \Id{d} + a$ and $G = \Id{d} + g$, one gets, simply, $dv = a + g$. So, for a given tensor $g$, the vector field $v$ must minimize an expression of the form
\begin{equation}
\label{eq:norm.growth}
\int_{\Omega} B\big(x, (dv + dv^T - g - g^T)/2\big) dx    
\end{equation}
where $B$ was discussed in \cref{sec:elastic.3D}. There is no loss of generality in assuming that $g$ is symmetric, which will be done in the following. The minimum of \cref{eq:norm.growth} is not always zero, i.e., the equation $\frac{dv + dv^T}{2} = g$ does not always have a solution. A necessary condition (which is sufficient when $\Omega$ is simply connected) is that $\nabla \times (g \times \nabla g) = 0$ (row-wise curl application, followed by column-wise; see, e.g., \citet{gonzalez2008first}).

\subsection{Riemannian viewpoint}
\label{sec:riemann.growth}
Returning to the Riemannian situation discussed in shape spaces, the metric was defined as
$\|v\|^2_\Omega = \kappa \|v\|_V^2 +  \llbracket v \rrbracket_{\Omega}^2$ with $\llbracket v \rrbracket_{\Omega}^2$ given by the right-hand side of \cref{eq:norm.growth} with $g=0$. One can apply the same approach here, letting
\[
\llbracket v \rrbracket_{\Omega}^2 = \inf_g \int_{\Omega} B\big(x, (dv + dv^T)/2 - g\big) dx.
\]
Obviously, this definition has little interest unless one restricts the space of growth tensors under consideration (otherwise, $\llbracket v \rrbracket_{\Omega}^2=0$ since one can take $g = (dv+dv^T)/2$). Letting $\mathcal G(\Omega)$ denote a set of tensor fields $(x \mapsto g(x) \in \Sym(\mR^d))$, one can define
\[
\llbracket v \rrbracket_{\Omega}^2 = \inf_{g\in \mathcal G(\Omega)} \int_{\Omega} B(x, (dv + dv^T)/2 - g) dx
\]
which is not trivial in general. If $\mathcal G (\Omega)$ forms a vector space, then $\llbracket v \rrbracket_{[\Omega]}$ is  a semi-norm on $V$.

Note that one can also switch the focus to the growth tensor and define, for $g\in \mathcal G(\Omega)$,
\[
\|g\|_\Omega^2 = \min_{v\in V} \left(\kappa \|v\|_V^2 + \int_{\Omega} B\big(x, (dv + dv^T)/2 - g\big) dx\right), 
\]
which defines a norm on growth tensors. The introduction of the regularization by the $V$ norm ensures that the minimum is attained at a unique $v\in V$, that one can denote $v_{g, \Omega}$, which depends linearly on $g$ and is such that $\kappa \|v_{g, \Omega}\|^2_V \leq \|g\|^2_{\Omega}$. One can therefore consider evolution equations in the form
\[
\left\{
\begin{aligned}
\partial_t \varphi(t, x) &= v_{g(t), \Omega(t)}(\varphi(t,x))\\
\Omega(t) &= \varphi(t, \Omega(0))
\end{aligned}
\right.
\]
which are well posed (starting with $\varphi(0, \cdot) = \id{3}$) as long as
\[
\int_0^1 \|g(t)\|_{\Omega(t)}^2 dt < \infty.
\]

This framework therefore provides two formally equivalent optimal control problems. In the first one, one minimizes, with respect to $v(\cdot)$
\begin{equation}
    \label{eq:riemann.growth.1}
\int_0^1 \|v(t)\|_{\Omega(t)}^2 dt,
\end{equation}
subject to $\varphi(1, \Omega_0) = \Omega_1$, $\varphi(0, \cdot) = \id{3}$, $\partial_t \varphi(t, \cdot) = v(t, \varphi(t, \cdot)$,  $\Omega(t) = \varphi(t, \Omega_0)$. In the second one, one minimizes, with respect to $g(\cdot)$
\begin{equation}
    \label{eq:riemann.growth.2}
\int_0^1 \|g(t)\|_{\Omega(t)}^2 dt,
\end{equation}
subject to $\varphi(1, \Omega_0) = \Omega_1$, $\varphi(0, \cdot) = \id{3}$, $\partial_t \varphi(t, \cdot) = v_{g(t)}(\varphi(t, \cdot))$, $g(t) \in \mathcal G(\Omega(t))$, $\Omega(t) = \varphi(t, \Omega_0)$. Both problems are, in addition, equivalent to minimizing, with respect to both $v(\cdot)$ and $g(\cdot)$, 
\begin{equation}
    \label{eq:riemann.growth.3}
\kappa \int_0^1 \|v(t)\|_{V}^2 dt + \int_0^1 \int_{\Omega(t)} B\big(x, (dv(t,x) + dv(t,x)^T)/2 - g(t,x)\big) dx
\end{equation}
%\textcolor{red}{shouldn't it be $\Omega(t)$ in the second integral (and in eq.(11) as well)?} 
subject to $\varphi(1, \Omega_0) = \Omega_1$, $\varphi(0, \cdot) = \id{3}$, $\partial_t \varphi(t, \cdot) = v(t, \varphi(t, \cdot))$, $g(t) \in \mathcal G(\Omega(t))$, $\Omega(t) = \varphi(t, \Omega_0)$.

When $\mathcal G(\Omega)$ is a vector space, the minimum value of these optimal control problems with given $\Omega_0$ and $\Omega_1$ is symmetric in $\Omega_0$ and $\Omega_1$ and its square root satisfies the triangular inequality. This minimum is always larger to that obtained with $B=0$ and therefore cannot be zero unless $\Omega_0 = \Omega_1$. (Note that the minimum can be infinite if the problem is unfeasible.)
Under suitable assumptions,  solutions of this optimal control problem always exist. A precise statement of this result and a sketch of its proof are provided in the appendix.

\subsection{Growth as an internal force}
Some additional notation is needed here.
Denote the topological dual of a Hilbert space $H$, with inner product $\scp{\cdot}{\cdot}_H$, by $H^*$ and if $\mu\in H^*$ is a linear form and if $h\in H$, denote their pairing by $\lform{\mu}{h}$ (i.e., $\mu(h)$). Riesz's representation theorem gives an isometric correspondence between $H$ and $H^*$ with, denoting by $K_H: H^* \to H$ the operator that associates to a linear form $\mu$ the unique vector $h\in H$ such that $\lform{\mu}{\tilde h} = \scp{h}{\tilde h}_H$ for all $\tilde h\in H$,  $\|h\|_H^2 = \lform{K_H^{-1}h}{h}$. This construction will be applied to $H=V$.

Introduce the (finite-dimensional) linear operator $\beta(x)$ operating on symmetric $3\times 3$ matrices such that $B(x, S) = \scp{S}{\beta(x) S}$ (with $\scp{S}{S'} = \trace(SS')$), and define, for a tensor field $x\mapsto S(x)$, 
\[
\boldsymbol \beta_\Omega(S) = \int_{\Omega} \beta(x) S(x) dx.
\]
Defining $\mathbbmss d v = (dv + dv^T)/2$, one has
\[
\int_{\Omega} B\big(x, (dv + dv^T)/2 - g\big) dx = \lform{\dd^*\boldsymbol \beta_\Omega\dd v}{v} - 2 \lform{\dd^* \boldsymbol\beta_{\Omega} g}{v} + \lform{\boldsymbol\beta_{\Omega} g}{g}.
\]
Letting $\mathbbm j_{g, \Omega} = \dd^* \boldsymbol\beta_\Omega g$, one has
\begin{equation}
    \label{eq:vg}
v_{g, \Omega} = (\kappa K_V^{-1} +  \dd^*\boldsymbol \beta_\Omega\dd)^{-1} \mathbbm j_{g, \Omega}.
\end{equation}

This relation provide an alternative way of modeling the growth process. One can indeed, following \citet{hsieh2022mechanistic}, directly define a ``yank'' (derivative of a force) $\mathbbm j$ as a control, with 
$v_{\mathbbm j} = (\kappa K_V^{-1} +  \dd^*\boldsymbol \beta_\Omega\dd)^{-1} \mathbbm j$ and use the running cost
\[
\int_0^1 \lform{\mathbbm j(t)}{v_{\mathbbm j(t)}} dt 
\]
with $\prt_t \varphi_{\mathbbm j}(t,x) =v_{\mathbbm j(t)}(\varphi_{\mathbbm j}(t,x))$. One can then show that the finiteness of the cost implies that the ODE has solutions over all time interval. One can also prove that optimal control $\mathbbm j$ always exist in this case. 

Note that this problem is different from the one described in equations \crefrange{eq:riemann.growth.1}{eq:riemann.growth.3}. In that case, one has
\[
\|g\|_\Omega^2 = \lform{\boldsymbol\beta_\Omega g}{g} - \lform{\mathbbm j_g}{v_{g}},
\]
showing that the  geodesics for the $\|\cdot\|_{\Omega}$ metric (which remain to be explored) are likely to behave differently than those studied in \citet{hsieh2022mechanistic}.

\subsection{A simple example}
\label{sec:example}
Assume that the growth tensor is scalar, i.e., $g(x) = \rho(x) \Id{3}$ and that $g(x) = 0$ on $\partial\Omega$, to avoid keeping track of boundary terms. Also assume that the elastic energy on $\Omega$ is homogeneous and isotropic (\cref{eq:isotropic}), which implies that $B(x, g(x))$ is proportional to $\rho(x)^2$, the proportionality constant being, using the Lam\'e coefficients, equal to $3(3\lambda/2 + \mu)$. Letting $\xi = 3\lambda/2+\mu$ and using the bilinearity of $B(x,\cdot)$ and the fact that $\trace(dv)=\trace(dv^T) = \nabla\cdot v$, a direct computation yields:
\[
 B\big(x, (dv + dv^T)/2 - g\big) = B(x, (dv + dv^T)/2) - 2 \xi \rho(x) \nabla\cdot v(x) + 3\xi \rho(x)^2
\]
Integrating by parts, one has
\[
\int_\Omega \rho(x) \nabla\cdot v(x) dx = - \int_\Omega \nabla \rho(x)^T v(x) dx
\]
so that, using the previous notation, 
\[
\mathbbm j_{g, \Omega} = -\xi \nabla \rho.
\]
One therefore finds that 
\[
\|\rho \Id{3}\|^2_\Omega = 3 \int_{\Omega} \rho(x)^2 dx - \int_{\Omega} \nabla\rho(x)^T (\kappa K_V^{-1} +  \dd^*\boldsymbol \beta_\Omega\dd)^{-1} \nabla\rho(x) dx.
\]

Similarly, the minimum in $\rho$ of $B\big(x, (dv(x) + dv(x)^T)/2 - \rho(x) \Id{3)}\big)$ is attained at $\rho = \nabla\cdot v/3$ and
\[
\|v\|^2_\Omega = \kappa \|v\|_V^2 +  \int_{\Omega} B\big(x, (dv(x) + dv(x)^T)/2 - \nabla\cdot v(x) /3\big) dx
\]
\subsection{Growth due to external action}

Shape variations resulting from a growth tensor as described above may be caused by external effects (e.g., impact of a disease) and do not need to follow a least-action principle such as described in the previous paragraph. More likely, the growth tensor will follow its own course, according to a process influenced by elements that are independent of the material properties of the deforming shape. The growth tensor evolution cannot be completely independent of the shape, however, since it must be supported by the time dependent domain $\Omega(t)$. It is also possible that changes in the geometry of the shape impact how growth behaves.

All this results in evolution systems with coupled evolution equations, typically involving moving domains. In \citet{bressan2018model}, a scalar growth is assumed, with the relationship $\nabla\cdot v = \rho$, consistent with \cref{sec:example}. The growth function depends on another function, $u$, representing the ``concentration of morphogen'', so that $\rho = \alpha\circ u$ for a fixed function $\alpha$. This morphogen concentration follows a partial differential equation (PDE), namely $\Delta u = w - u $, with  Neumann's boundary conditions, where $w$ itself is a density advected by the motion, i.e., satisfying $\prt_t w + \nabla\cdot(vw) = 0$, which provides the coupling between growth and shape change. Initial conditions are the initial domain $\Omega_0$ and the initial value of $w$, $w_0$. One can then show that, when starting with a domain $\Omega_0$ with smooth enough boundary and with a smooth enough density $w_0$, a solution to the growth system exists over some finite interval $[0, T]$ for some (small enough) $T$. 

In \citet{hsieh2022mechanistic} and \citet{hsieh2021diffeomorphic}, the additional regularization term $\|v\|_V$ described in this chapter is added, using the formulation in \cref{eq:vg}
\[
(\kappa K_V^{-1} +  \dd^*\boldsymbol \beta_\Omega\dd) v = \mathbbm j,
\]
where $\mathbbm j$ is the modeled control
(as seen in our simple example of \cref{sec:example}, $\mathbbm j$ has an interpretation similar to that of $-\nabla\rho$). \citet{hsieh2022mechanistic} models $\mathbbm j$ as a function $\mathbbm j(\varphi, \theta)$, for some time-independent parameter $\theta$, providing coupled equations 
\[
\left\{
\begin{aligned}
&\prt_t \varphi(t,x) = v(t, \phi(t,x))\\
& (\kappa K_V^{-1} +  \dd^*\boldsymbol \beta_\Omega\dd) v(t, \cdot) = \mathbbm j(\phi(t, \cdot), \theta)
\end{aligned}
\right.
\]
The system is shown to have a unique solution $t \mapsto \phi(t, \cdot)$ over arbitrary large time intervals, for any fixed $\theta$, provided that $\mathbbm j(\phi, \theta)$ is Lipschitz in $\phi$ for the $(1, \infty)$ norm. Denoting this solution by $\phi(t, \cdot;\theta)$, this property allows for the specification of optimization problems over the parameter $\theta$ involving the transformation $\phi(1, \cdot;\theta)$. 
\medskip

\begin{figure}
    \centering
    \includegraphics[trim=4cm 3cm 4cm 2cm, clip, width=0.32\textwidth]{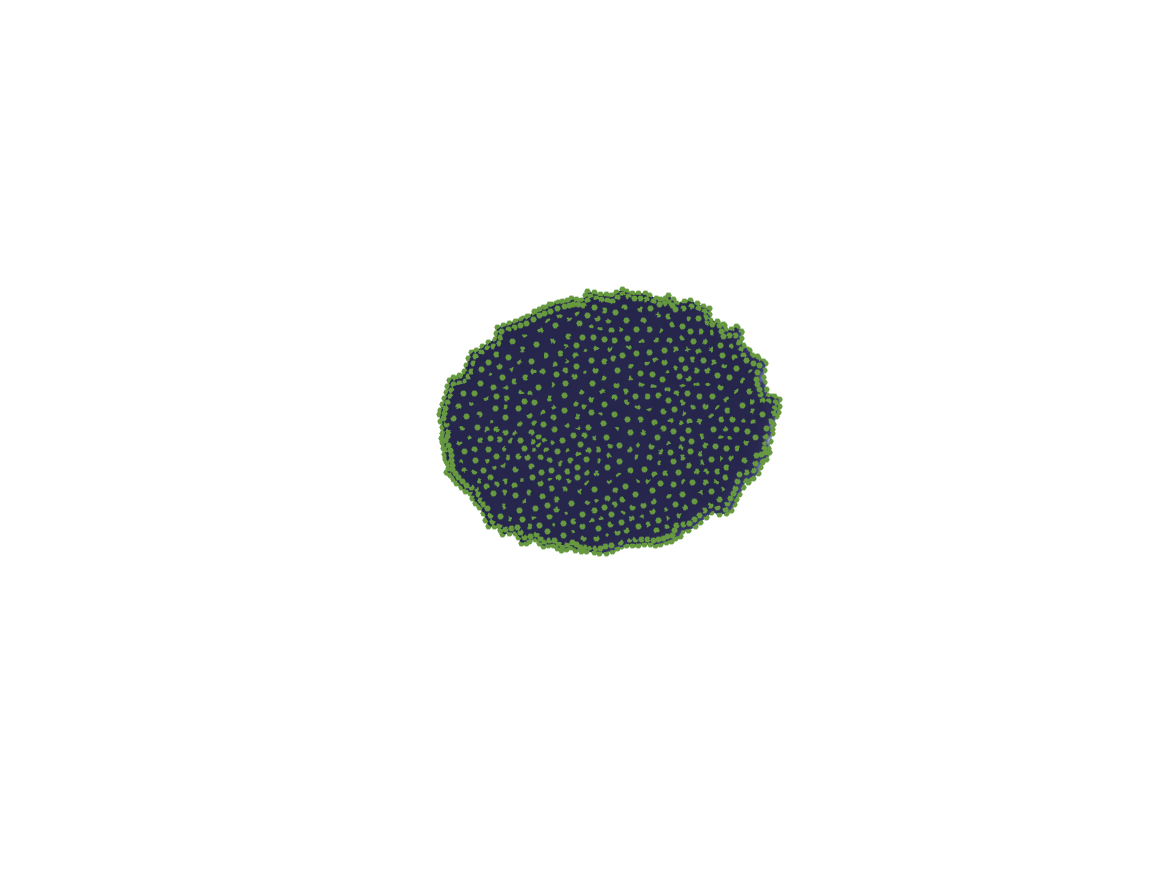}
    \includegraphics[trim=4cm 3cm 4cm 2cm, clip, width=0.32\textwidth]{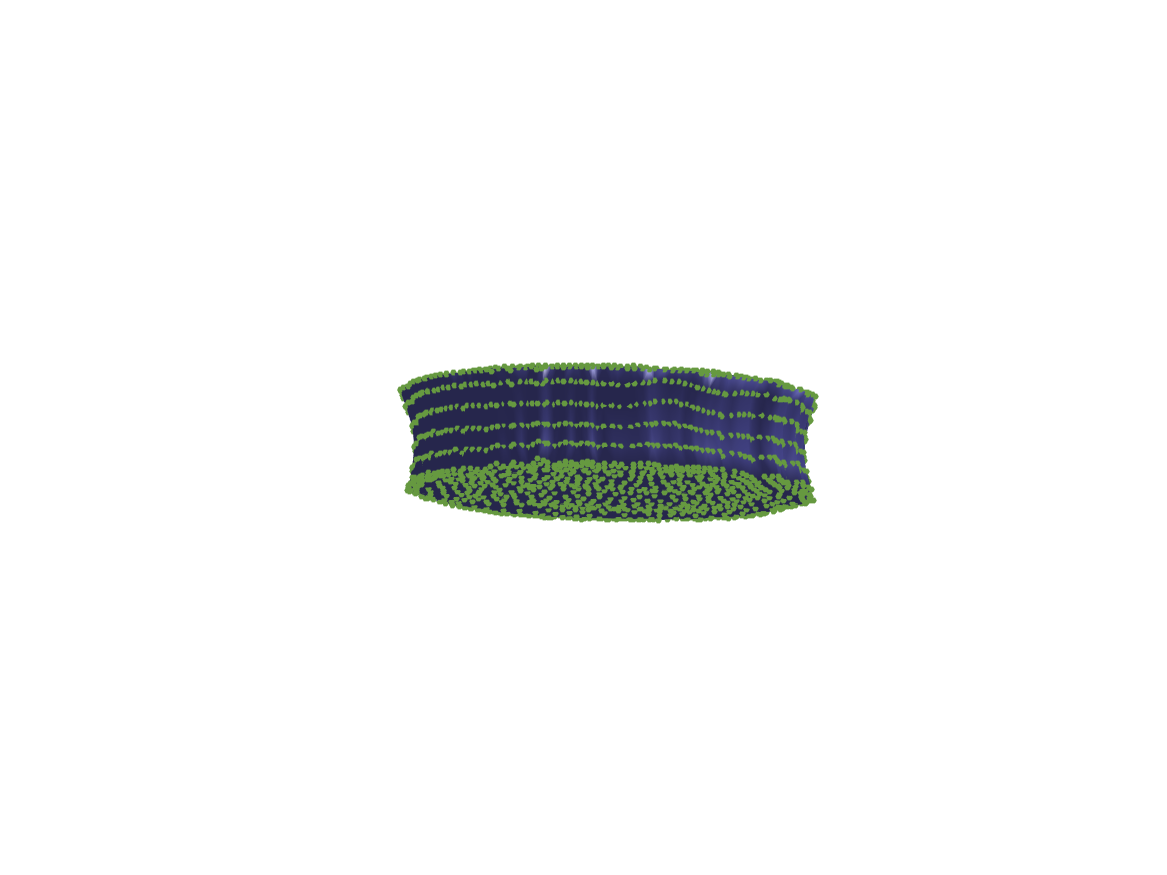}\\
    \includegraphics[trim=4cm 3cm 4cm 2cm, clip, width=0.32\textwidth]{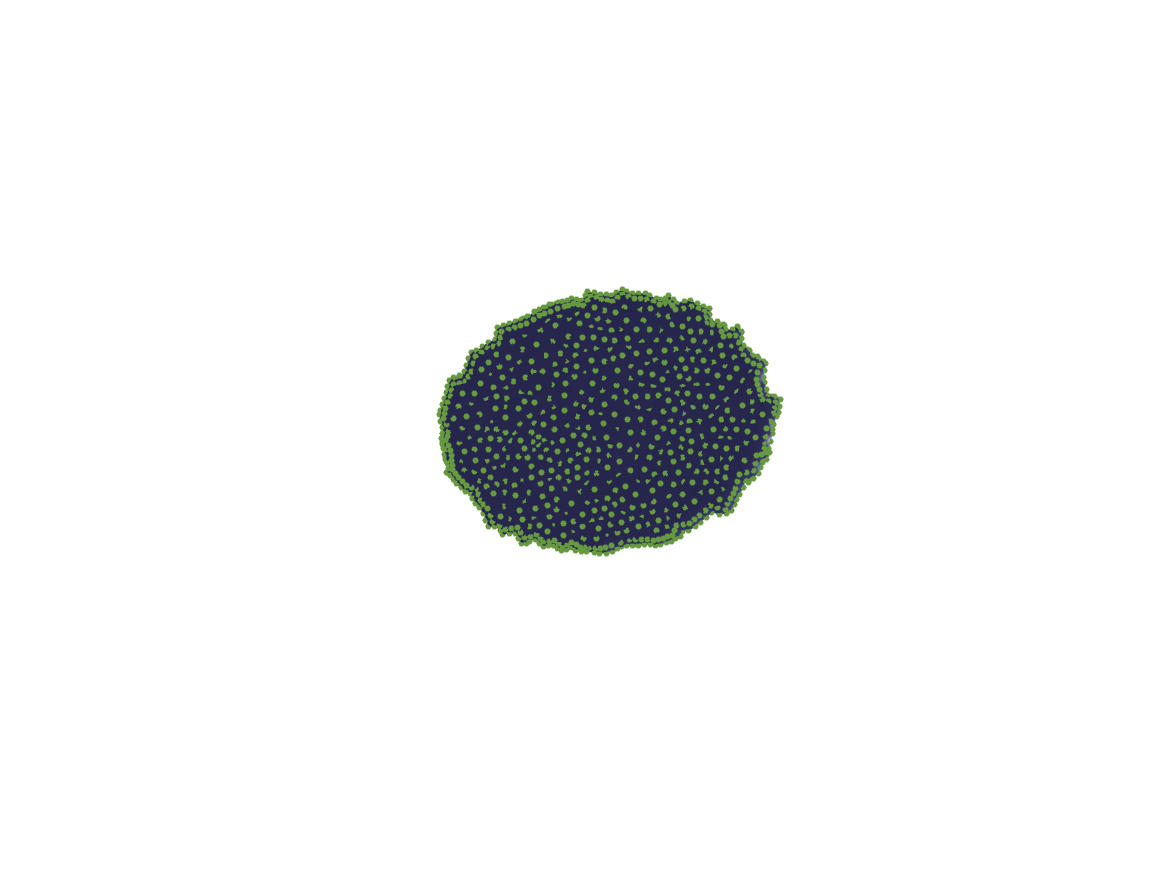}
    \includegraphics[trim=4cm 3cm 4cm 2cm, clip, width=0.32\textwidth]{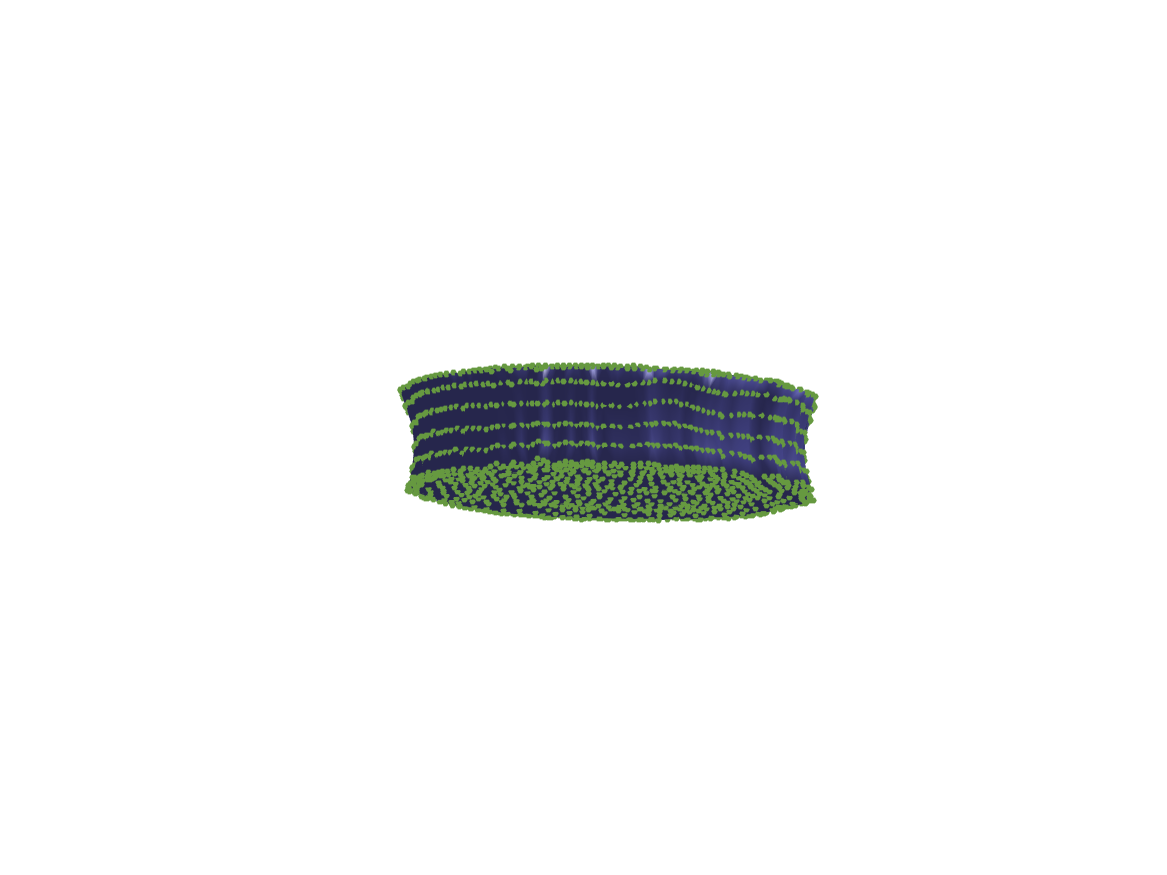}\\
    \includegraphics[trim=4cm 3cm 4cm 2cm, clip, width=0.32\textwidth]{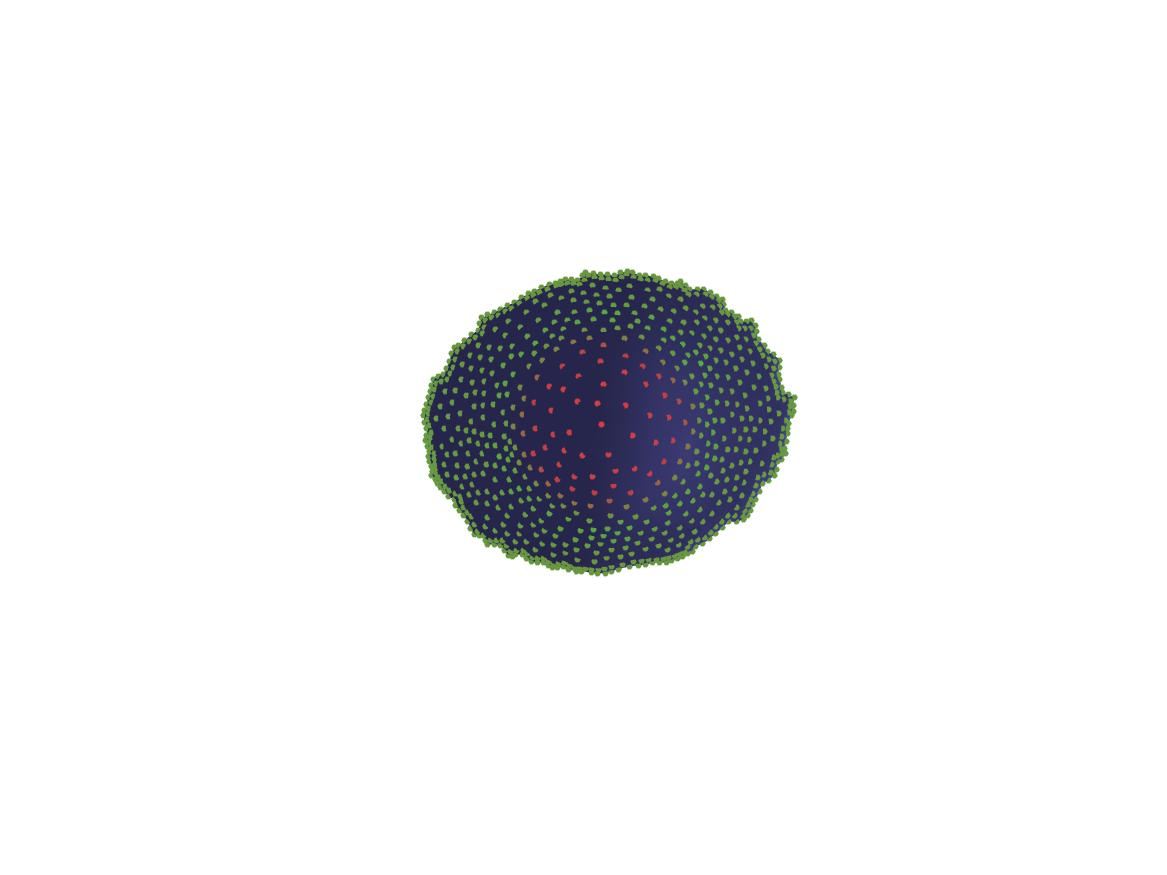}
    \includegraphics[trim=4cm 3cm 4cm 2cm, clip, width=0.32\textwidth]{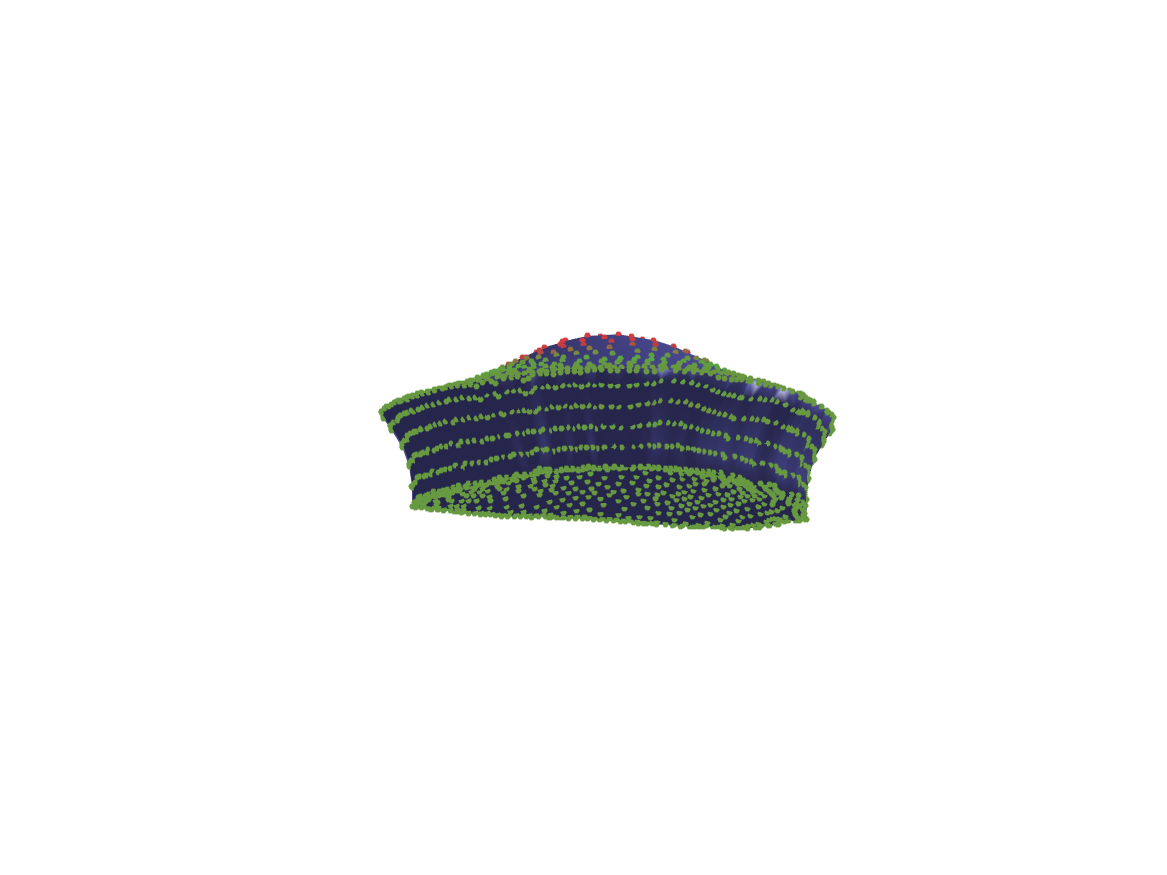}\\
    \includegraphics[trim=4cm 3cm 4cm 2cm, clip, width=0.32\textwidth]{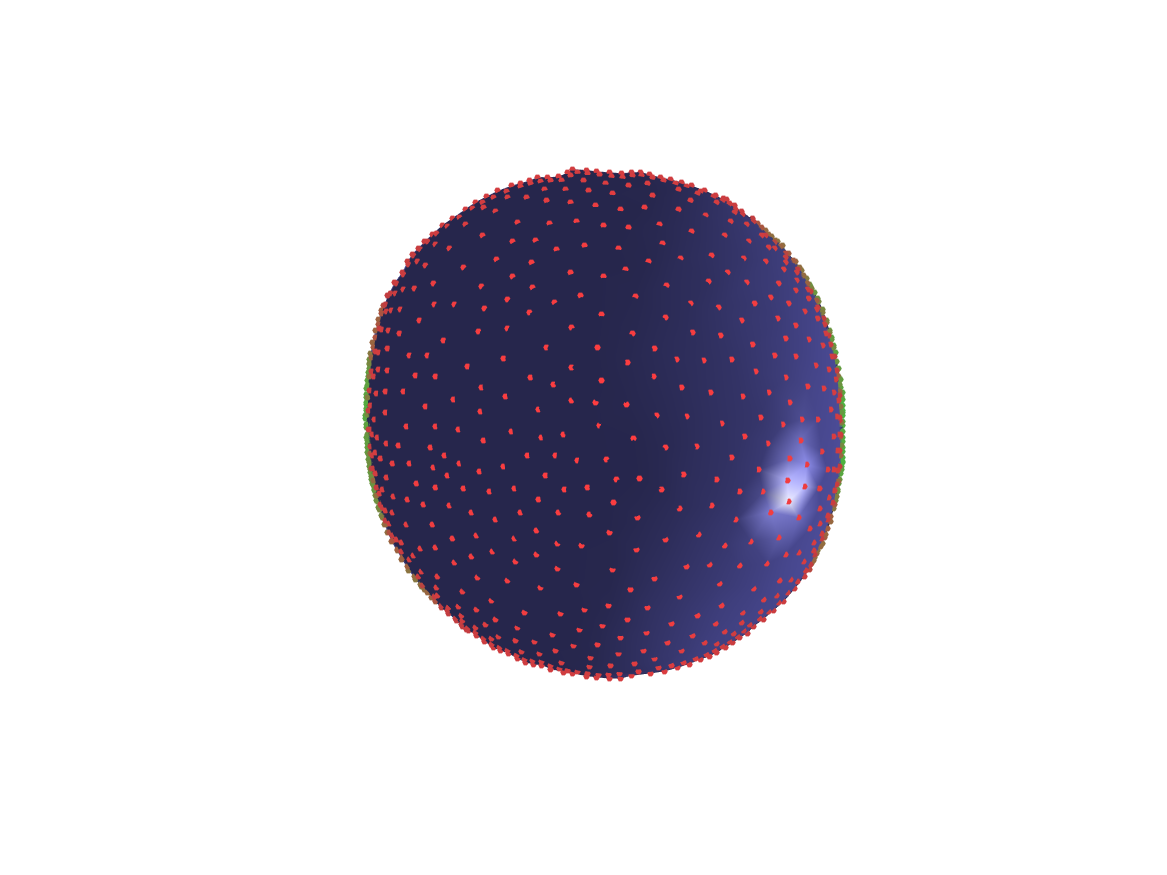}
    \includegraphics[trim=4cm 3cm 4cm 2cm, clip, width=0.32\textwidth]{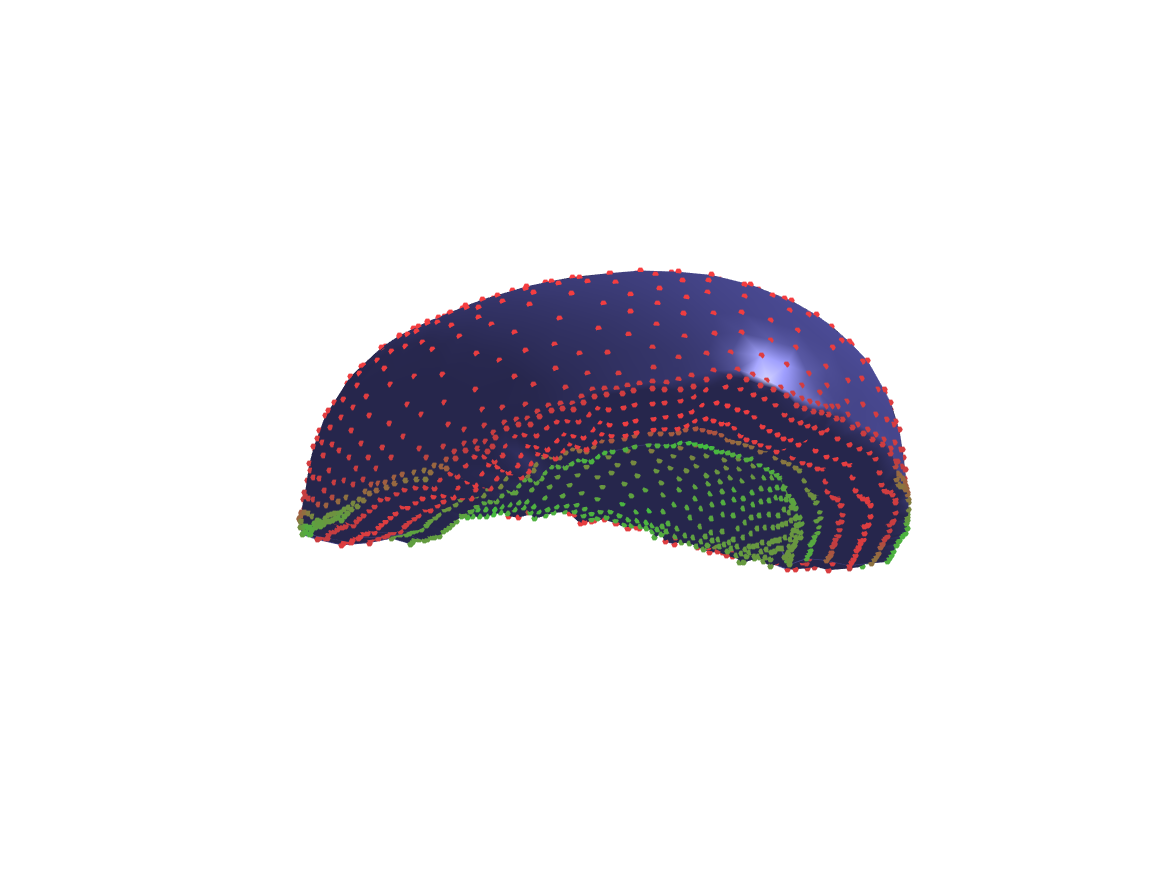}
    \caption{Growth model from \citet{hsieh2021model-based} applied to a 3D volume. Dots are colored proportionally to the magnitude of $p$ in \cref{eq:growth.hsieh}. Rows 1 to 4 provide two views of the evolving shapes at times $t=0$, $t=0.33$, $t=0.67$ and $t=1.0$. (Images generated from code developed by Dai-Ni Hsieh.)}
    \label{fig:growth.hsieh}
\end{figure}

A more complex system is introduced in \citet{hsieh2021diffeomorphic} in which $\mathbbm j$ is itself modeled based on a solution of a ``reaction-diffusion-convection'' equation on the moving domain $\Omega$. Ignoring a few technicalities, $\mathbbm j$ is given by $\mathbbm j = \nabla(Q( p))$ where $Q$ is a fixed function and $ p$ satisfies
\[
\partial_t p = \nabla\cdot(S_\varphi\nabla p - pv) + R(p)
\]
where $R$, the reaction function, is fixed, and $S_\varphi$, the diffusion matrix, is allowed to evolve with the transformation $\varphi$. One can then formulate suitable conditions under which the system
\begin{equation}
\label{eq:growth.hsieh}
\left\{
\begin{aligned}
&\prt_t \varphi(t,x) = v(t, \phi(t,x))\\
& (\kappa K_V^{-1} +  \dd^*\boldsymbol \beta_\Omega\dd) v(t, \cdot) = \nabla(Q( p))\\
&\partial_t p = \nabla\cdot(S_\varphi\nabla p - pv) + R(p)
\end{aligned}
\right.
\end{equation}
has solutions over arbitrary time intervals for a given initialization $p_0 = p(0, \dots)$. The determination of this initial condition for an optimal behavior at time 1 is tackled in  \citet{hsieh2021model-based}, where the existence of solutions of the optimization problem is shown. \Cref{fig:growth.hsieh} provides an example of growth process obtained as solution of this system.

\subsection{Constraints, deformation modules and other growth models}
Specific behavior can be enforced in a deformation process by constraining the values of the vector field at given locations in the shape. Theoretical bases for constrained and sub-Riemannian versions of LDDMM were introduced in \citet{arguillere_shape_2014,arguillere2015shape,arguillere2017sub-riemannian}, and a survey of such methods is provided in \citet{younes_sub-riemannian_2020}. Among such  approaches, deformation modules  \citep{gris2018sub-riemannian,lacroix2021imodal} offer a generic framework in which various types of behaviors can be defined by combining suitable constraints in a modular manner. Referring to the publications above for more details, the example of ``implicit elastic modules'' is closely related to this chapter's discussion. For such modules, the vector field $v$ is obtained as a minimizer of
\[
v \mapsto \lambda \|v\|_V^2 + \sum_{k=1}^m |\epsilon_v(x_k) - S_k(h_k)|^2
\]
where $\epsilon_v = (dv + dv^T)/2$, $x_1, \ldots, x_N$ are control points that are attached to (and move together with) the evolving shape, and $h \mapsto S_k(h)$ are symmetric matrices, parametrized by a control $h$, inducing a desired behavior (e.g., dilation) near the control points. This norm therefore introduces a finite set of (soft) constraints on the strain tensor.

A different approach at modeling growth can be found in \citet{kaltenmark2016geometrical,kaltenmark2019estimation}. In this work, a growing shape at a given time $t$ is defined as a transformation $q_t$ of a co-dimension-one foliation $X$, which encodes the full growth process. During the evolution, only the restriction of $q_t$ to the set $X_t$ formed by leaves at time $s\leq t$ of the foliation is relevant to describe the growing shape. The value of $q_t(x)$ remains constant until $t$ is reaches the foliation index of $x$, so that the function $q_0$ encodes all future initializations of the growth process. This process can be constructed through an evolution equation in the form $\partial_t q_t = v(t, q_t)$, and an example is developed in \citet{kaltenmark2019estimation} to model animal horn growth.

\section{Conclusion}
Starting from the notion of shape spaces built along the principles of Grenander's metric pattern theory and the action of diffeomorphism groups, this chapter surveyed a few recent efforts to incorporate physical constraints in the modelling of trajectories in such spaces. It first discussed the class of hybrid models that consist in combining the original shape space metric induced by the deformation group with other more physically-informed metrics, in particular those derived from linear elasticity theory. A second general approach is to further constrain shape evolution via the introduction of a growth model underlying the morphological transformation. 

One of the main motivation behind all of these works is to advance the ability of shape space frameworks to model physical or biological processes, while still preserving the advantages of the geometric shape space metric setting. Indeed, this enables the formulation of the dynamics of those processes as control systems and provides adequate regularization norms to ensure existence and smoothness of solutions in many cases. Furthermore, by considering the associated optimal control problems, those same models can often lead to natural and well-posed approaches for tackling the inverse problem of e.g. determining the causes/sources of morphological changes based on some observed shape evolution. The ideas provided by the present chapter are examples of emerging efforts toward cross-fertilization between the fields of shape analysis, mathematical biology, biomedical engineering and material science.        

\bibliography{references}

\appendix
\section{Elastic surface metric as the limit of the laminar model (\cref{sec:surfaces})}
\label{app:laminar}
\begin{figure}
    \centering
    \includegraphics[width=0.5\textwidth]{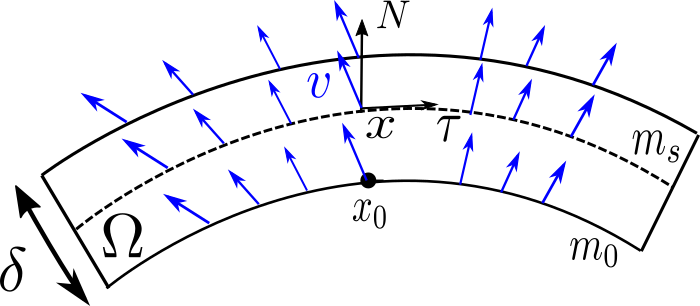}
    \caption{Cross-sectional schematic representation of the thin-shell layered elastic domain with the deformation field $v$ in blue.}
    \label{fig:2D.elastic.domain}
\end{figure}

Given an oriented surface $\mathcal M_0$ in $\mR^3$, and its unit normal vector field denoted $\nu_0$, one can generate a foliated 3D volume as the set of points $\varPhi(s,x_0) = x_0 + s \delta \nu_0(x_0)$, $x_0\in \mathcal M_0$, $s\in [0,1]$, and $\varPhi$ is a diffeomorphism for small enough $\delta > 0$. In this case, the unit normal $N$ to the layer $\mathcal M_s = \varPhi(\{s\} \times M_0)$ at the point $x = \varPhi(x_0,s) \in \Omega$ is also $N(x) = \nu_0(x_0)$. It coincides, up to a factor $\delta$, with $S = \prt_s \varPhi$ and satisfies $dN N = 0$. Let $v_0: \mathcal M_0 \to \mR^3$ be a vector field on $\mathcal M_0$, and define its extension $v$ to $\Omega$ by $v(\varPhi(s,x_0)) = v_0(x_0)$, so that $v$ satisfies $dv N = 0$. See  \cref{fig:2D.elastic.domain} for an illustration. Let $\sigma_0 = d\nu_0$ denote the shape operator on the surface $\mathcal M_0$ and similarly $\sigma_s$ the shape operator of layer $\mathcal Ms$ i.e. the restriction of $dN$ to the tangent space of $\mathcal M_s$. Recall that the shape operator on a surface is a symmetric operator.  

Write $v = v_T + v_N N$ where $v_T \in \mR^3$ is tangent to the layers and $v_N$ is scalar. If $\tau$, $\tilde\tau$ are vectors tangent to the layers, we have
\begin{equation}
\label{eq:computation_elastic_energy.surf.inter0}
\tilde\tau^T dv \tau = \tilde \tau^T dv_T \tau + (\nabla v_N^T\tau) (\tilde\tau^T N) + v_N \tilde\tau^T dN\tau = 
\tilde \tau^T dv_T \tau + v_N \tilde\tau^T dN\tau
\end{equation}
In particular, letting $\varepsilon_T = (dv_T + dv_T^T)/2$, and for $\{\tau_1, \tau_2\}$ an orthonormal basis of the tangent plane to $\mathcal M_s$ at $x$, one has:
\[
\tau_1^T \varepsilon \tau_1 + \tau_2^T \varepsilon \tau_2 = \tau_1^T dv_T \tau_1 + \tau_2^T dv_T \tau_2 + v_N [\tau_1^T dN \tau_1 + \tau_2^T dN \tau_2].
\]
With the sum of the first two terms, one recognizes the divergence of $v_T$ on the surface $\mathcal M_s$  which will be denoted by $\nabla_{\mathcal M_s} \cdot v_T$. Similarly, the term within brackets is the divergence of the shape operator on $\mathcal M_s$ which equals $-2 H_{\mathcal M_s}$ where $H_{\mathcal M_s}$ is the mean curvature of $\mathcal M_s$. Therefore, one deduces that, on $\mathcal M_s$:
\[
\tau_1^T \varepsilon \tau_1 + \tau_2^T \varepsilon \tau_2 = \nabla_{\mathcal M_s} \cdot v_T - 2v_N H_{\mathcal M_s}.
\]

Moreover, as $dv N = 0$, it follows that $N^T \varepsilon N = 0$ and thus, on $\mathcal M_s$:
\begin{equation*}
    \trace(\varepsilon) = \tau_1^T \varepsilon \tau_1 + \tau_2^T \varepsilon \tau_2 + N^T \varepsilon N = \nabla_{\mathcal M_s} \cdot v_T - 2v_N H_{\mathcal M_s}.
\end{equation*}

Similarly, looking at the second term in \cref{eq:computation_elastic_energy} and using \cref{eq:computation_elastic_energy.surf.inter0}, one has:
\begin{align*}
&(\tau_1^T \varepsilon \tau_1)^2 + (\tau_2^T \varepsilon \tau_2)^2 + 2(\tau_1^T \varepsilon \tau_2)^2 \\
&= (\tau_1^T dv_T \tau_1 + v_N \tau_1^T \sigma_s \tau_1)^2 + (\tau_2^T dv_T \tau_2 + v_N \tau_2^T \sigma_s \tau_2)^2 + 2 (\tau_1^T \varepsilon_T \tau_2 + v_N \tau_1^T \sigma_s \tau_2)^2\\
&=(\tau_1^T dv_T \tau_1)^2 + (\tau_2^T dv_T \tau_2)^2 + 2 (\tau_1^T \varepsilon_T \tau_2)^2 \\
&\phantom{=}+v_N^2\left[(\tau_1^T \sigma_s \tau_1)^2 + (\tau_2^T \sigma_s \tau_2)^2 +2 (\tau_1^T \sigma_s \tau_2)^2 \right]\\
&\phantom{=}+2v_N\left[(\tau_1^T dv_T \tau_1)(\tau_1^T \sigma_s \tau_1) + (\tau_2^T dv_T \tau_2)(\tau_2^T \sigma_s \tau_2) + 2 (\tau_1^T \varepsilon_T \tau_2) (\tau_1^T \sigma_s \tau_2) \right].
\end{align*}
In this computation, one uses the fact that the operator $dN$ restricted to the to the tangent space to $\mathcal M_s$ at $x$ (i.e. the space spanned by $\tau_1$ and $\tau_2$) coincides with $\sigma_s$. Now, by symmetry, one has $\tau_1^T dv_T \tau_1 = \tau_1^T \varepsilon_T \tau_1$ and $\tau_2^T dv_T \tau_2 = \tau_2^T \varepsilon_T \tau_2$. Moreover, recalling that for any $2\times 2$ symmetric tensors $\omega$ and $\tilde{\omega}$, one has $\trace(\omega \tilde{\omega}) = \omega_{1,1} \tilde{\omega}_{1,1} + \omega_{2,2} \tilde{\omega}_{2,2} +2\omega_{1,2} \tilde{\omega}_{1,2}$, one gets:
\begin{align*}
 (\tau_1^T \varepsilon \tau_1)^2 + (\tau_2^T \varepsilon \tau_2)^2 + 2(\tau_1^T \varepsilon \tau_2)^2&=\trace(\varepsilon_T^2) + v_N^2\trace(\sigma_s^2) + 2 v_N \trace(\varepsilon_T \sigma_s) \\
 &=\trace((\varepsilon_T+v_N\sigma_s)^2).
\end{align*}

Now, using the symmetry of $\varepsilon$ and the fact that $N^T \varepsilon N = 0$:
\begin{equation}
\label{eq:computation_elastic_energy.surf.inter1}
N^T \varepsilon^2 N = |\varepsilon N|^2 = (\tau_1^T \varepsilon N)^2 + (\tau_2^T \varepsilon N)^2.
\end{equation}
If $\tau$ is tangent to the layers, one has $\tau^T dv N = 0$ and 
\begin{align}
\label{eq:computation_elastic_energy.surf.inter2}
\tau^T dv^T N = N^T dv \tau &= N^T dv_T \tau + (\nabla v_N^T\tau) (N^T N) + v_N N^T dN\tau \nonumber\\
&= 
N^T dv_T \tau + \nabla v_N^T\tau .
\end{align}
Moreover, since $N^T v_T = 0$, it follows that $N^T dv_T \tau = - v_T^T dN \tau$. Using this together with \cref{eq:computation_elastic_energy.surf.inter1}, \cref{eq:computation_elastic_energy.surf.inter2}, and with the fact that $dN$ is symmetric, one deduces that: 
\begin{align*}
    N^T \varepsilon^2 N - (N^T \varepsilon N)^2&= ((-dN^T v_T + \nabla v_N)^T\tau_1)^2 + ((-dN^T v_T + \nabla v_N)^T\tau_2)^2 \\
    &= |-\sigma_s v_T + \nabla_{\mathcal M_s} v_N|^2
\end{align*}
where $\nabla_{\mathcal M_s}$ is the gradient operator on $\mathcal M_s$. 

Based on all the above expressions, one can finally rewrite \cref{eq:computation_elastic_energy} at $x = \varPhi(x_0, s)$ as
\begin{align}
	\label{eq:computation_elastic_energy.surf}
	B(x, \varepsilon)=&\lambda_{\mathrm{tan}} \left(\nabla_{\mathcal M_s} \cdot v_T - 2H_{\mathcal M_s} v_N\right)^2 
		+\mu_{\mathrm{tan}}\trace((\varepsilon_T + v_N \sigma_s)^2)
		\\ \nonumber
		&	+ 2 \, \mu_{\mathrm{ang}} \left|-\sigma_s v_T + \nabla_{\mathcal M_s} v_N\right|^2 ,
\end{align}
and using by a change of variables in the integral expression of the energy, one further has:
\[
\frac1\delta \int_{\Omega} B(x, \varepsilon) dx = \frac1\delta \int_0^{1} \int_{\mathcal M_0} B(x_0 + s \delta \nu_0, \varepsilon) |J_{\varPhi}(s,x_0)| d\mathrm{vol}_{m_0}(x_0) ds
\]
where $|J_{\varPhi}(s,x_0)|$ denotes the Jacobian determinant of $\varPhi$ at $(s,x_0)$. As $\partial_s \varPhi(s,x_0) = \delta \nu_0(x_0)$ and $d_{x_0} \varPhi(s,x_0) = \text{Id} + s \delta d\nu_0(x_0)$, one gets $d_{x_0} \varPhi(0,x_0) = \text{Id}$ where $\text{Id}$ denotes here the identity on the tangent space to $m_0$ at $x_0$. Therefore, $|J_{\varPhi}(0,x_0)|=\delta$ for all $x_0 \in m_0$. Consequently, taking the limit $\delta\to 0$ in the above and using the continuity of $B$ and $J_{\varPhi}$ leads to the following expression of the elastic metric on the surface $\mathcal M_0$: 
\[
\llbracket v \rrbracket_{\mathcal M_0}^2 = \int_{m_0} B(x_0, \varepsilon) d\mathrm{vol}_{\mathcal m_0}(x_0)
\]
with $B$ given by \cref{eq:computation_elastic_energy.surf} (with $s=0$). Furthermore, it can be easily checked, based on their expressions in the frame $(\tau_1,\tau_2,N)$, that the three terms in $B(x_0,\varepsilon)$ correspond precisely, up to multiplicative constants, to the ones of \cref{eq:surf.elastic.metric} thus showing that the elastic metric in \cref{eq:surf.elastic} can be also recovered as the thin shell limit of the 3D laminar model introduced in  \cref{sec:elastic.3D}. 

\section{Existence of optimal paths (\cref{sec:riemann.growth})}
Considering the minimization problem introduced in \crefrange{eq:riemann.growth.1}{eq:riemann.growth.3}, this section proves that, under suitable assumptions, optimal solutions exist. These  assumptions are as follows.
\begin{enumerate}[label=(\arabic*),wide]
    \item Let $p\geq 1$. The Hilbert space $V$ is continuously embedded in the Banach space $C^p_0(\mR^3, \mR^3)$ of $p$ times continuously differentiable vector fields that vanish (with their first $p$ derivatives) at infinity, with the norm 
    \[
    \|v\|_{p, \infty} = \sum_{k=0}^p \max\{|d^kv(x)|:, x\in \mR^3\}.
    \]
    \item $V$ is also continuously embedded in $H^1(\mR^3, \mR^3)$, the Sobolev space of square-integrable functions with square-integrable first derivatives. 
    \item The mapping $x \mapsto B(x, \cdot)$ from $\mR^3$ to the set of positive semi-definite quadratic forms is continuous in $x$. In particular, $|B(x, \cdot)|$ is bounded on compact subsets of $\mR^3$. 
    \item There exists a constant $c$ such that 
    $B(x, S) \geq c |S|^2$ for all $S\in \Sym$ and all $x \in \mR^3$.
            \label[condition]{item:app.2}
    \item The sets $\mathcal G(\Omega)$, defined over  compact subsets of $\Omega \subset \mR^3$, satisfy the following conditions.
    \begin{enumerate}[label=\theenumi-\roman*,leftmargin=0.5cm,itemindent=1cm]
    \item If $\Omega\subset \tilde\Omega$, then $\mathcal G(\Omega) \subset \mathcal G(\tilde \Omega)$.
    \item Define, for $\delta>0$, $\Omega^\delta = \{x: \mathrm{dist}(x, \Omega) \leq \delta\}$. Then $\bigcap_{\delta > 0} \mathcal G(\Omega^\delta) = \mathcal G(\Omega)$.
    \label[condition]{item:app.3.ii}
        \item $\mathcal G(\Omega)$ is a strongly closed convex subset of $H_\Sym := L^2(\mR^3, \Sym(\mR^3))$.
            \label[condition]{item:app.3.iii}
% \item $\mathcal G(\Omega) = \Phi_\Omega(\Theta)$ from a fixed compact parameter set $\Theta\subset \mR^k$, with the continuity property that if $\psi_n$ is a sequence of diffeomorphisms that tend to $\id{3}$ for the $1, \infty$ norm, then $\Phi_{\psi_n(\Omega}(\theta)$ converges to $\Phi_{\Omega}(\theta)$ for all $\theta\in \Theta$.
    \end{enumerate}
    \label[condition]{item:app.3}
    % \item For every $\varphi$ diffeomorphism of $\mR^3$, there exists a one-to-one mapping $\tau_\varphi: x \mapsto \tau_\varphi(x)$ defined on $\mR^3$ and taking values in the set of linear transformations of $\Sym$ such that, for all $\Omega$ compact in $\mR^3$, 
    % \[
    % \mathcal G(\varphi(\Omega)) = \{\tau_\phi(x) g \circ \varphi^{-1}\},
    % \]
    % with $\|\tau_\varphi - \mathrm{Id}_{\Sym}\|_\infty \leq \|\varphi - \id{3}\|_{1, \infty}$.
    % \label[condition]{item:app.4}
\end{enumerate}
For example the sets $\mathcal G(\Omega) = \{g\, \Id{3}: g\in L^2(\Omega)\}$ satisfy \cref{item:app.3}.

% A first consequence of these assumptions is that, for any $v\in V$, the stress tensor $dv + dv^T$ is continuous and tends to zero at infinity. The mapping 
% \[
% g \mapsto \boldsymbol B_{\Omega} (dv+dv^T - 2g) := \int_{\Omega} B(x, dv(x) + dv(x)^T - 2g(x)) dx
% \]
% is convex and continuous in $L^2(\mR^3, \Sym)$ and therefore admits a minimizer under  \cref{item:app.3}. 

\medskip
Making these assumptions, let $v_n(\cdot)\in L^2([0,1], V)$ and $g_n\in L^2([0,1], H_\Sym)$ be minimizing sequences for the considered problem. To shorten notation, let $\varepsilon_n = (dv_n + dv_n^T)/2$.
Because $v_n$ is bounded in $L^2([0,1], V)$, one can replace it by a subsequence that converges weakly to some $v$ in that space, and using arguments developed in \citet{dupuis1998variation,trouve1995action,younes2019shapes}, the flows $\varphi_n$ associated with $v_n$ converge uniformly in time and uniformly on compact sets in space to the flow $\varphi$ associated with $v$. From weak convergence and weak lower semicontinuity of the norm, one has
\[
\int_0^1 \|v\|_V^2 dt \leq \liminf \int_0^1 \|v_n\|_V^2 dt
\]
and from the convergence of the flows, one has $\varphi(1, \Omega_0) = \Omega_1$ because this holds for each $\varphi_n$. 
%Indeed, since $\varphi_n(1, x) \in \Omega_1$ for all $x\in \Omega_0$, the fact that $\varphi(1, \Omega_0) \subset \Omega_1$ results from the closedness of $\Omega_1$. Conversely, if $y\in \Omega_1$, there exists a sequence $x_n\in \Omega_0$ such that $\varphi_n(1, x_n) = y$. Since $\Omega_0$ is compact, one can replace $x_n$ by a converging subsequence to some $x\in \Omega_0$ and $\varphi_n(1, x_n)$ converges to $\varphi(1,x)$. \textcolor{red}{Do we really need to justify this? I think it is pretty clear given the uniform convergence of the $\varphi_n$ no?}

% We note that, for all $x\in \mR^3$ and $S\in \Sym$,
% \[
% c |S|^2 \leq B(x, S) \leq \left(B(x, \varepsilon - S)^{1/2} + B(x, \varepsilon)^{1/2}\right)^2 \leq 2B(x, \varepsilon - S) + 2B(x, \varepsilon)
% \]
% Moreover, for all $n>0$
% \[
% \int_0^1 \int_{\Omega_n(t)}B(x, \varepsilon) dx dt \leq C \int_0^1 |v|_{H^1}^2 dt  \leq C' \int_0^1 |v|_{V}^2 dt
% \]
% for some constant $C,C'$ that do not depend on $n$. This implies that
% \[
% \int_0^1 \|g_n\|_2^3 dt
% \]
% is bounded and that we can assume, using a subsequence if needed that $g_n \wto g$ in $L^1([0,1], H_\Sym)$.

Based on the assumptions made on $B$, one has, for all $x \in \mR^3$ and $t\in[0,1]$:
\begin{align*}
c|g_n(t,x)|^2 \leq B(x,g_n(t,x)) &\leq \left(B(x,\varepsilon_n(t,x) -g_n(t,x))^{1/2} + B(x,\varepsilon_n(t,x))^{1/2} \right)^2\\
&\leq 2\left(B(x,\varepsilon_n(t,x) -g_n(t,x)) + B(x,\varepsilon_n(t,x)) \right)
\end{align*}

Because of the convergence of $\phi_n$, there exists a compact set $\bar\Omega \subset \mR^d$ that contains all the $\Omega_n(t)$, $n \in \mathbb{N}$, $t\in [0,1]$. This implies that there exist constants $C,C'$ such that, for all $n \in \mathbb{N}$ (using the boundedness of $B(x,\cdot)$ on compact sets):
\[
\int_0^1 \int_{\Omega_n(t)} B(x,\varepsilon_n(t,x)) dx dt \leq C \int_0^1 \int_{\Omega_n(t)} |\varepsilon_n(t,x)|^2 dx dt \leq C' \int_0^1 \|v_n(t)\|_{H^1}^2 dt.
\]
By the continuous embedding of $V$ into $H^1$, $\|v_n(t)\|_{H^1}$ is bounded up to a multiplicative constant by $\|v_n(t)\|_V$, which implies that the above term is bounded independently of $n$. The same holds for:
\[
\int_0^1 \int_{\Omega_n(t)} B(x,\varepsilon_n -g_n) dx dt = \frac{1}{4} \int_0^1 \int_{\Omega_n(t)} B(x,dv_n +dv_n^T -2g_n) dx dt
\]
as $(v_n,g_n)$ is a minimizing sequence for the functional in \cref{eq:riemann.growth.3}. This implies that the sequence $\int_0^1 \|g_n\|^2_{H_\Sym} dt$ is bounded and that one can assume, using a subsequence if needed, that $g_n \wto g$ in $L^2([0,1], H_\Sym)$.

It remains to prove that $g(t) \in \mathcal G(\Omega(t))$ to show that $(v, g)$ provides a solution of the minimization problem.
Fixing $\delta >0$, one can restrict the minimizing sequence to those large enough $n$ for which $\max\{|\phi_n(t,x) - \phi(t,x)|, t\in [0,1], x\in \bar\Omega\} < \delta$, so that $\Omega_n(t) \subset \Omega^\delta(t)$ for all $n$ and $t$.

Let
\[
\Gamma(\Omega(\cdot), \delta) = \{\tilde g(\cdot): \tilde g(t) \in \mathcal G(\Omega^\delta(t)), \ \text{for a.e } t\in [0,1]\},
\]
so that $g_n \in \Gamma(\Omega(\cdot), \delta)$.
This is a convex set, which follows directly from our hypotheses on the sets $\mathcal G(\Omega)$, and it is closed in $L^2([0,1], H_\Sym)$. Indeed, if $\tilde g_n \in \Gamma(\Omega(\cdot), \delta)$ converges to $\tilde g \in L^2([0,1], H_\Sym)$,  then a subsequence converges for almost all $t\in[0,1]$ and since each $\mathcal G(\Omega^\delta(t))$ is closed in $H_\Sym$, it results that $\tilde g(t) \in \mathcal G(\Omega^\delta(t))$ for almost all $t$. Now, as strongly closed convex sets are also weakly closed in $L^2([0,1], H_\Sym)$ (see, \citet{hytonen2016analysis}), one deduces from $g_n \wto g$ that $g \in \Gamma(\Omega(\cdot), \delta)$. Since this is true for all $\delta >0$, one has, taking a sequence $\delta_n\to 0$, that $g(t) \in \mathcal G(\Omega(t))$ for almost all $t\in [0,1]$.

This concludes the proof that $(v,g)$ is a minimizer of \cref{eq:riemann.growth.3}.

\end{document}